\documentclass[3pt]{article}

\usepackage{graphicx}\usepackage{epstopdf}\usepackage{color}
\usepackage{amsmath}\usepackage{amsfonts}\usepackage{amssymb}
\numberwithin{equation}{section}

\begin{document}

\title{Two-grid discretizations and a local finite element scheme for a non-selfadjoint Stekloff eigenvalue problem}
\author{{Hai Bi, Yu Zhang, Yidu Yang*}\\
\small{School of Mathematical Science, Guizhou Normal University,}\\
\small{GuiYang, 550001, China}\\
\small{bihaimath@gznu.edu.cn, zhang\_hello\_hi@126.com, ydyang@gznu.edu.cn}}
\date{~} \pagestyle{plain}
\maketitle

\indent{\bf\small Abstract~:~} {\small In this paper, for a new Stekloff eigenvalue problem which is non-selfadjoint and not $H^1$-elliptic, we establish and analyze two kinds of two-grid discretization scheme and a local finite element scheme. We present the error estimates of approximations of two-grid discretizations. We also prove a local error estimate which is suitable for the case that the local
refined region contains singular points lying on the boundary of domain. Numerical experiments are reported finally to show the efficiency of our schemes.
}\\
\indent{\bf\small Keywords~:~\scriptsize} {\small
Stekloff eigenvalue problems, two-grid discretizations, local error
estimates, local computation.}

\section{Introduction}
\label{intro} \indent  Numerical methods for solving Stekloff
eigenvalue problems have attracted the attention of academia for their important physical background
and wide applications. Till now,
systematical and profound studies on the finite element approximations are mainly for selfadjoint Stekloff eigenvalue
problems, for example, see \cite{aa1,andreev,armentano2004,armentano2008,armentano2011,bramble,cao,cheng,garau2011,huang,lqx,lim,russo,xie}, etc.\\
 \indent Recently the study for Stekloff eigenvalues in inverse scattering has aroused researchers' interest (see \cite{cakoni2016,liusun}).
The differential operator corresponding to this problem is non-selfadjoint and the associated weak formulation does not satisfy $H^1$-elliptic condition, which are the main differences from those studied before.
\cite{cakoni2016} studies the mathematical properties of this problem and its conforming finite element approximation, later \cite{liusun} proves the error estimate of eigenvalues. In this paper we study further the finite element method for the problem, and the features of our work are as follows.\\
\indent (1) The existing work analyzes that the discrete solution operator $T_h$ converges to $T$, the solution operator of source problem, in $L^{2}(\partial\Omega)$. In this paper, to devise and analyze efficient schemes, we prove that $T_h$ converges to $T$ in "negative" space $H^{-\frac{1}{2}}(\partial\Omega)$, then we give the error estimates of eigenfunctions and eigenvalues. We also give the local a priori error estimates. With the local error estimates we establish and analyze a local computational scheme.\\
\indent (2) The two-grid discretization introduced by Xu \cite{xu1992,xu1996}  is an efficient method for reducing the computational costs and maintaining the accuracy of approximation at the same time. This powerful computing technique has been used and developed by many scholars later (see, e.g., \cite{cai2009,chien2006,hu2011,mu2007,xu2,yang2011,zhzsc}). In this paper we establish two kinds of two-grid discretization scheme for the Stekloff eigenvalue problem, in particular, the second scheme performs better because the matrices are constructed to be symmetric and definite in solving linear systems. We provide the error analysis and numerical experiments to show the efficiency of our schemes. \\
\indent (3) For elliptic boundary value problems, Xu and Zhou \cite{xu1} combine the two-grid
finite element discretizations with the local defect-correction technique to
propose a local and parallel finite element algorithm, and this computing technique has been applied successfully to
other problems (see, e.g., \cite{bi,dai,he2010,he2008}). For the eigenfunctions with local low smoothness, or singularity, based on the two-grid discretizations we present a local computational scheme. Theoretical analysis and numerical experiment all indicate that the local correction does work as we expected.\\
\indent The rest of this paper is organized as follows. In the next
section, some preliminary materials are presented. Local a priori error estimates for conforming finite elements approximations for the Stekloff eigenvalue problem are analyzed in Section 3. In Section 4, two kinds of two-grid discretization scheme
for the Stekloff eigenvalue problem are established and analyzed. In Section 5, a local finite element scheme is presented and
its error estimates are proved. Numerical experiments are provided in Section 6 to show the efficiency of our schemes.\\
\indent We refer to  \cite{babuska,boffi,brenner,ciarlet,oden,sun}  as regards the basic theory of finite element methods in this paper.\\
\indent Throughout this paper, $C$ denotes a generic positive constant independent of mesh diameters, which
may not be the same at each occurrence. For simplicity, we use the symbol $a\lesssim b$ to mean that
$a\leq C b$.

\section{Preliminaries}
\label{sec:1}
\indent Let $\Omega \subset \mathbb{R}^{2}$ be a bounded polygon with Lipshitz boundary $\partial\Omega$ and
$\nu$ be the unit outward normal to $\partial\Omega$. Let $H^{t}(\Omega)$ denote the usual Sobolev
space with real order $t$ on $\Omega$  and $H^{0}(\Omega)=L_{2}(\Omega)$. $\|\cdot\|_{t,\Omega}$ is the norm
on $H^{t}(\Omega)$. Let $H^{t}(\partial\Omega)$ denote the Sobolev space with real order
$t$ on $\partial\Omega$ with the norm $\|\cdot\|_{t,\partial\Omega}$.\\
\indent Consider the following Stekloff eigenvalue problem to find $\lambda\in \mathbb{C}$ and a nontrivial function
$u\in H^{1}(\Omega)$ such that
\begin{eqnarray}\label{s2.1}
&& \left \{
\begin{array}{ll}
\triangle u +k^{2}n(x) u=0&~~~in~ \Omega,\\
\frac{\partial u}{\partial \nu}+\lambda u=0&~~~on~ \partial \Omega,
\end{array}
\right.
\end{eqnarray}
where $k$ is the wavenumber and $n(x)$ is the index of refraction.
Assume that  $n=n(x)$ is a bounded complex valued function given by
$$n(x)=n_{1}(x)+i\frac{n_{2}(x)}{k},$$
where $i=\sqrt{-1}$, $n_{1}(x)>0$ and $n_{2}(x)\geq 0$ are bounded smooth functions.\\
Denote
$$(u,v)=\int\limits_{\Omega} u \overline{v}dx,~~~\langle f,g\rangle=\int\limits_{\partial\Omega} f \overline{g}ds,$$
and define the continuous sesquilinear form
$$a(u,v):=(\nabla u,\nabla v)-k^{2}(nu,v),~~~\forall u,v\in H^{1}(\Omega).$$
For any $g\in  H^{1}(\Omega)$, $\langle f,g\rangle$ has a continuous extension to $f\in H^{-\frac{1}{2}}(\partial\Omega)$ so that
$\langle f,g\rangle$ is continuous on $H^{-\frac{1}{2}}(\partial\Omega)\times H^{\frac{1}{2}}(\partial\Omega)$.\\

\indent The weak form of (\ref{s2.1}) is to find $(\lambda,u)\in
\mathbb{C}\times H^{1}(\Omega)$, $u\neq 0$, such that
\begin{eqnarray}\label{s2.2}
a(u,v)=-\lambda \langle u,v\rangle,~~~\forall v\in H^{1}(\Omega).
\end{eqnarray}
From \cite{liusun} we know that $a(\cdot,\cdot)$ satisfies G{\aa}rding's inequality, i.e., there exist constants
$K<\infty$ and $\alpha_{0}>0$ such that
\begin{eqnarray*}
Re\{a(v,v)\}+K\|v\|_{0,\Omega}^{2}\geq \alpha_{0}\|v\|_{1,\Omega}^{2},~~~\forall v\in H^{1}(\Omega).
\end{eqnarray*}
Let $K$ be a positive constant which is large enough, and define the sesquilinear form
\begin{eqnarray*}
\widetilde{a}(u,v):=a(u,v)+K(u,v)=(\nabla u,\nabla v)-k^{2}(nu,v)+K(u,v),~u,v\in H^{1}(\Omega),
\end{eqnarray*}
then it is easy to verify that $\widetilde{a}(\cdot,\cdot)$ is $H^{1}(\Omega)$-elliptic (see \cite{liusun}).\\

\indent Let $\pi_{h}=\{\tau\}$ be a mesh of $\Omega$, $h(x)$ be the diameter
of the element $\tau$ containing $x$, $h_{G}=\max\limits_{x\in
G}h(x)$, and $h_{\Omega}=h$ be the mesh diameter of $\pi_{h}$. Let $V_{h}(\Omega)\subset C(\overline{\Omega})$, defined on $\pi_{h}$, be a piecewise polynomial space of degree $m(m\geq 1)$ and $V_{h}^{B}:=V_{h}(\Omega)|_{\partial\Omega}$ be the restriction of $V_{h}(\Omega)$ on $\partial\Omega$.\\
\indent We assume that the finite element spaces in this paper satisfy the following
condition (see, e.g.,  \cite{xu1}): \\
\indent (C0)~{\em Approximation.} There exists $m\geq 1$ such that for $\psi\in H^{1+s}(\Omega)$,
\begin{eqnarray*}
\inf\limits_{v\in V_{h}(\Omega)}(\|h^{-1}(\psi-v)\|_{0,\Omega}+\|\psi-v\|_{1,\Omega})\leq
Ch^{s}\|\psi\|_{1+s,\Omega},~~~0\leq s\leq m.
\end{eqnarray*}

\indent The finite element approximation of (\ref{s2.2}) is to find $(\lambda_{h},u_{h})\in \mathbb{C}\times V_{h}(\Omega)$, $u_{h}\neq 0$ such that
\begin{eqnarray}\label{s2.5}
a(u_{h},v_{h})=-\lambda_{h} \langle u_{h},v_{h}\rangle,~~~\forall v_{h}\in V_{h}(\Omega).
\end{eqnarray}

\indent Consider the following source problem (\ref{s2.6}) associated with
(\ref{s2.2}), and the approximate source problem (\ref{s2.7}) associated with (\ref{s2.5}):\\
Given $g\in H^{-\frac{1}{2}}(\partial\Omega)$, find $\phi\in H^{1}(\Omega)$ such that
\begin{eqnarray}\label{s2.6}
a(\phi,v)=\langle g,v\rangle,~~~\forall v \in H^{1}(\Omega);
\end{eqnarray}
Find $\phi_{h} \in V_{h}(\Omega)$ such that
\begin{eqnarray}\label{s2.7}
a(\phi_{h},v)=\langle g,v\rangle,~~~\forall v \in V_{h}(\Omega).
\end{eqnarray}
Introduce the following Neumann eigenvalue problem:
\begin{eqnarray}\label{s2.8}
&& \left \{
\begin{array}{ll}
\triangle \phi +k^{2}n(x) \phi=0&~~~in~ \Omega,\\
\frac{\partial \phi}{\partial \nu}=0&~~~on~ \partial \Omega.
\end{array}
\right.
\end{eqnarray}
When $k^2$ is not a Neumann eigenvalue of (\ref{s2.8}), from Fredholm Alternative (see, e.g., Section 5.3 of \cite{gc}) we can prove that for $g\in H^{-\frac{1}{2}}(\partial\Omega)$, there exists
a unique solution $\phi\in H^1(\Omega)$ to (\ref{s2.6}) satisfying
\begin{eqnarray}\label{s2.9}
\|\phi\|_{1,\Omega}\leq C\|g\|_{-\frac{1}{2},\partial\Omega}.
\end{eqnarray}
Thus, one can define the operator $A: H^{-\frac{1}{2}}(\partial\Omega)\to H^{1}(\Omega)$ by
\begin{eqnarray}\label{s2.10}
a(Ag,v)=\langle g,v\rangle,~~~\forall v\in H^{1}(\Omega),
\end{eqnarray}
and the Neumann-to-Dirichlet map $T: H^{-\frac{1}{2}}(\partial\Omega)\rightarrow H^{\frac{1}{2}}(\partial\Omega)$ by
\begin{eqnarray}\label{s2.11}
Tg=(Ag)',
\end{eqnarray}
where $'$ denotes the restriction to $\partial\Omega$, namely, $Tg=Ag|_{\partial\Omega}$.\\
Then, (\ref{s2.2}) has the equivalent operator form as follows:
\begin{eqnarray}\label{s2.12}
Au=-\frac{1}{\lambda}u.
\end{eqnarray}
Similarly, (\ref{s2.7}) defines a discrete operator $A_{h}:H^{-\frac{1}{2}}(\partial\Omega)\to V_{h}(\Omega)$ satisfying
\begin{eqnarray*}
a(A_{h}g,v)=\langle g,v\rangle,~~~ \forall v \in V_{h},
\end{eqnarray*}
and $T_{h}: H^{-\frac{1}{2}}(\partial\Omega)\rightarrow V_h^B$ such that
\begin{eqnarray}\label{s2.13}
T_{h}g=(A_{h}g)'.
\end{eqnarray}
Then, (\ref{s2.5}) has the equivalent operator form as follows:
\begin{eqnarray}\label{s2.14}
A_{h}u_{h}=-\frac{1}{\lambda_{h}}u_{h}.
\end{eqnarray}
\indent In this paper, we always assume that $k^2$ is not a Neumann eigenvalue of (\ref{s2.8}).\\
\indent Consider the dual problem of (\ref{s2.2}): Find $(\lambda^*,u^*)\in
\mathbb{C}\times H^{1}(\Omega)$, $u^*\neq 0$, such that
\begin{eqnarray}\label{s2.15}
a(v,u^*)=-\overline{\lambda^*}\langle v,u^*\rangle,~~~\forall v\in H^{1}(\Omega).
\end{eqnarray}
The primal and dual eigenvalues are connected via $\lambda=\overline{\lambda^*}$.\\
\indent The discrete variational formulation associated with (\ref{s2.15}) is to find $(\lambda_{h}^*,u_{h}^*)\in \mathbb{C}\times V_{h}(\Omega)$, $u_{h}^*\neq 0$ such that
\begin{eqnarray}\label{s2.16}
a(v_{h}, u_{h}^*)=-\overline{\lambda_{h}^*}\langle v_{h}, u_{h}^*\rangle,~~~\forall v_{h}\in V_{h}(\Omega).
\end{eqnarray}
The primal and dual eigenvalues are connected via $\lambda_{h}=\overline{\lambda_{h}^*}$.\\
\indent Similarly, from the source problem corresponding to (\ref{s2.15}) and (\ref{s2.16}) we can define the operator $A^*$, $T^*$ and $A_{h}^*$, $T_{h}^*$, respectively.\\

\indent For (\ref{s2.6}), there holds the following regular estimates which will be used in the sequel.\\
\indent{\bf Lemma 2.1.} Let $\phi$ be the
solution of (\ref{s2.6}). If $g\in L^{2}(\partial\Omega)$, then
$\phi\in H^{1+\frac{r}{2}}(\Omega)$ and
\begin{eqnarray}\label{s2.17}
\|\phi\|_{1+\frac{r}{2},\Omega}\leq C\|g\|_{0,\partial\Omega};
\end{eqnarray}
If $g\in H^{\frac{1}{2}}(\partial\Omega)$, then $\phi\in
H^{1+r}(\Omega)$ and
\begin{eqnarray}\label{s2.18}
\|\phi\|_{1+r,\Omega}\leq
C\|g\|_{\frac{1}{2},\partial\Omega}.
\end{eqnarray}
Here $r=1$ when the largest inner angle $\theta$ of $\Omega$ satisfying $\theta<\pi$, and $r<\frac{\pi}{\theta}$ which can
be arbitrarily close to $\frac{\pi}{\theta}$ when $\theta>\pi$.\\
\indent{\bf Proof.} According to \cite{dauge1988} and Proposition 4.1 in \cite{aa1} we can prove the desired results.~~~$\Box$\\
\indent Lemma 2.1 guarantees that the eigenfunction of (\ref{s2.2}) $u\in H^{1+r}(\Omega)$. \\

Let $P_{h}: H^{1}(\Omega)\to
 V_{h}(\Omega)$ be the projection defined by
\begin{eqnarray}\label{s2.19}
a(\phi-P_{h}\phi, v)=0,~~~\forall v\in V_{h}(\Omega).
\end{eqnarray}
Thus, for any $g \in H^{-\frac{1}{2}}(\Omega)$, we have
\begin{eqnarray*}
&&a(A_{h}g-P_{h}(Ag),v)=a(A_{h}g-Ag+Ag-P_{h}(Ag),v)\\
&&~~~~~~=a(A_{h}g-Ag,v)+a(Ag-P_{h}(Ag),v)=0,~~~ \forall v\in V_{h}(\Omega).
\end{eqnarray*}
Since the above equation admits a unique solution, we have $A_{h}g=P_{h}Ag, \forall g \in H^{-\frac{1}{2}}(\Omega)$, then $A_{h}=P_{h}A$.\\

\indent{\bf Lemma 2.2.} Let $\phi$ be the solution of
(\ref{s2.6}). If $\phi\in
H^{1+t}(\Omega)$~$(t\geq r)$, then
\begin{eqnarray}\label{s2.20}
&~&\|\phi-P_{h}\phi\|_{1,\Omega}\lesssim h^{\sigma}\|\phi\|_{1+t,\Omega},\\\label{s2.21}
&~&\|\phi-P_{h}\phi\|_{-\frac{1}{2},\partial\Omega}\lesssim h^{r+\sigma}\|\phi\|_{1+t,\Omega},
\end{eqnarray}
where $\sigma=\min\{m,t\}$ and the principle to determine $r$ see Lemma 2.1.\\
{\bf{Proof.}} From Theorem 3.1 in \cite{liusun} and the interpolation error estimates we can immediately get (\ref{s2.20}).\\
Next we use Aubin-Nitsche technique to prove (\ref{s2.21}). According to the definition of $A$ we deduce that for any $v\in V_{h}(\Omega)$,
\begin{eqnarray*}
&&\|\phi-P_{h}\phi\|_{-\frac{1}{2},\partial\Omega}=\sup\limits_{g\in H^{\frac{1}{2}}(\partial\Omega),g\neq 0}\frac{|\langle \phi-P_{h}\phi,g\rangle|}{\|g\|_{\frac{1}{2},\partial\Omega}}\\
&&~~~=\sup\limits_{g\in H^{\frac{1}{2}}(\partial\Omega),g\neq 0}\frac{|\overline{\langle g, \phi-P_{h}\phi\rangle}|}{\|g\|_{\frac{1}{2},\partial\Omega}}\\
&&~~~=\sup\limits_{g\in H^{\frac{1}{2}}(\partial\Omega),g\neq 0}\frac{|a(Ag,\phi-P_{h}\phi)|}{\|g\|_{\frac{1}{2},\partial\Omega}}\\
&&~~~=\sup\limits_{g\in H^{\frac{1}{2}}(\partial\Omega),g\neq 0}\frac{|\overline{a(\phi-P_{h}\phi,Ag)}|}{\|g\|_{\frac{1}{2},\partial\Omega}}\\
&&~~~=\sup\limits_{g\in H^{\frac{1}{2}}(\partial\Omega),v\in V_{h}(\Omega)}\frac{|a(\phi-P_{h}\phi,Ag-v)|}{\|g\|_{\frac{1}{2},\partial\Omega}}\\
&&~~~\leq\sup\limits_{g\in H^{\frac{1}{2}}(\partial\Omega),v\in V_{h}(\Omega)}\frac{\|\phi-P_{h}\phi\|_{1,\Omega}\|Ag-v\|_{1,\Omega}}{\|g\|_{\frac{1}{2},\partial\Omega}}\\
&&~~~\lesssim h^{r}\|\phi-P_{h}\phi\|_{1,\Omega},
\end{eqnarray*}
where the last inequality is valid because of the interpolation estimate and (\ref{s2.18}). Substituting  (\ref{s2.20}) into the above inequality we obtain (\ref{s2.21}). The proof is completed. ~~~$\Box$\\

Since we also need the error estimate $\|\phi-P_{h}\phi\|_{0,\Omega}$, now we consider an auxiliary problem:
Find $\xi_{f} \in H^{1}(\Omega)$, such that
\begin{eqnarray}\label{s2.33}
a(v,\xi_{f})=(v,f),~~~\forall v\in H^{1}(\Omega).
\end{eqnarray}
\indent From Theorem 2.1 and Remark 3 in \cite{babuska2} we have the following regularity result.\\
\indent{\bf Lemma 2.3.} If $f\in L^{2}(\Omega)$, then there exists a
unique solution $\xi_{f}\in H^{1+r}(\Omega)$ to (\ref{s2.33}), and
\begin{eqnarray}\label{s2.34}
\|\xi_{f}\|_{1+r,\Omega}\leq C\|f\|_{0,\Omega},
\end{eqnarray}
where the principle to determine $r$ see Lemma 2.1.\\

\indent{\bf Lemma 2.4.} For any $\phi\in H^{1}(\Omega)$, there holds
\begin{eqnarray}\label{s2.35}
\|\phi-P_{h}\phi\|_{0,\Omega} & \lesssim & h^{r}\|\phi\|_{1,\Omega},\\\label{s2.36}
\|\phi-P_{h}\phi\|_{0,\Omega} & \lesssim & h^{r}\|\phi-P_{h}\phi\|_{1,\Omega}.
\end{eqnarray}
{\bf{Proof.}} Let $\widetilde{P_{h}}: L^2(\Omega)\rightarrow V_{h}(\Omega)$ be the projection defined by $$a(v, \xi-\widetilde{P_{h}}\xi)=0, \forall v\in V_{h}(\Omega).$$
\noindent Taking
$$f=\phi-P_{h}\phi ~~and ~~v=\phi-P_{h}\phi$$
in (\ref{s2.33}), and by using Aubin-Nitsche's technique we get
\begin{eqnarray*}
&&\|\phi-P_{h}\phi\|_{0,\Omega}^2=(\phi-P_{h}\phi, \phi-P_{h}\phi)=a(\phi-P_{h}\phi,\xi_{\phi-P_{h}\phi})\\
&&~~~~~~=a(\phi-P_{h}\phi,\xi_{\phi-P_{h}\phi}-\widetilde{P_{h}}\xi_{\phi-P_{h}\phi})\\
&&~~~~~~=a(\phi, \xi_{\phi-P_{h}\phi}-\widetilde{P_{h}}\xi_{\phi-P_{h}\phi})\\
&&~~~~~~\lesssim \|\phi\|_{1,\Omega}\|\xi_{\phi-P_{h}\phi}-\widetilde{P_{h}}\xi_{\phi-P_{h}\phi}\|_{1,\Omega}\\
&&~~~~~~\lesssim \|\phi\|_{1,\Omega}\cdot h^{r}\|\xi_{\phi-P_{h}\phi}\|_{1+r,\Omega}\\
&&~~~~~~\lesssim \|\phi\|_{1,\Omega}\cdot h^{r}\|\phi-P_{h}\phi\|_{0,\Omega},
\end{eqnarray*}
i.e.,
\begin{eqnarray*}
\|\phi-P_{h}\phi\|_{0,\Omega}\lesssim h^{r}\|\phi\|_{1,\Omega}.
\end{eqnarray*}
From the above deduction we can easily get (\ref{s2.36}). This ends the proof.~~~$\Box$\\

\indent From (\ref{s2.35}) we can get the following property of the projection $P_{h}$ which is obvious in the case that $a(\cdot,\cdot)$ is coercive. But, unfortunately  $a(\cdot,\cdot)$ in this paper is not coercive. \\
\indent{\bf Lemma 2.5.} For any $\phi\in H^{1}(\Omega)$, there holds
\begin{eqnarray}\label{s2.37}
\|P_{h}\phi\|_{1,\Omega} \lesssim \|\phi\|_{1,\Omega}.
\end{eqnarray}
{\bf{Proof.}} Since
\begin{eqnarray*}
a(\phi-P_{h}\phi, P_{h}\phi)=0,
\end{eqnarray*}
we derive that
\begin{eqnarray*}
|a(P_{h}\phi, P_{h}\phi)|=|a(\phi, P_{h}\phi)|\lesssim \|\phi\|_{1,\Omega}\|P_{h}\phi\|_{1,\Omega},
\end{eqnarray*}
thus,
\begin{eqnarray*}
&&\|P_{h}\phi\|_{1,\Omega}^{2}\lesssim |a(P_{h}\phi, P_{h}\phi)+K(P_{h}\phi, P_{h}\phi)|\\
&&~~~~~~\lesssim \|\phi\|_{1,\Omega}\|P_{h}\phi\|_{1,\Omega}+K\|P_{h}\phi\|_{0,\Omega}\|P_{h}\phi\|_{1,\Omega},
\end{eqnarray*}
so we have
\begin{eqnarray*}
&&\|P_{h}\phi\|_{1,\Omega}\lesssim \|\phi\|_{1,\Omega}+\|P_{h}\phi\|_{0,\Omega}.
\end{eqnarray*}
Noting that $\|P_{h}\phi\|_{0,\Omega}\leq \|P_{h}\phi-\phi\|_{0,\Omega}+\|\phi\|_{0,\Omega}$, we just need to prove $\|P_{h}\phi-\phi\|_{0,\Omega}\lesssim  \|\phi\|_{1,\Omega}$, which follows from (\ref{s2.35}).~~~$\Box$ \\

\indent (\ref{s2.9}) can be expressed as $\|Ag\|_{1,\Omega}\leq C\|g\|_{-\frac{1}{2},\partial\Omega}$, thus, from $A_{h}=P_{h}A$ and (\ref{s2.37}) we have
$\|A_{h}g\|_{1,\Omega}\leq C\|g\|_{-\frac{1}{2},\partial\Omega}$.\\

\indent With the error estimates of boundary value problem, (\ref{s2.20}), (\ref{s2.21}) and (\ref{s2.36}), we can get the error estimates of eigenvalue problem according to the classical spectral approximation theory (see \cite{babuska}) as long as we prove that $\|T-T_{h}\|_{H^{-\frac{1}{2}}(\partial\Omega)\rightarrow H^{-\frac{1}{2}}(\partial\Omega)}\rightarrow 0$. \\
\indent{\bf Lemma 2.6.} $\|T-T_{h}\|_{H^{-\frac{1}{2}}(\partial\Omega)\rightarrow H^{-\frac{1}{2}}(\partial\Omega)}\rightarrow 0$ as $h\rightarrow 0$ and $T$ is a compact operator.\\
{\bf{Proof.}} From the definitions of $A$, $A_{h}$, $T$ and $T_{h}$ we have
\begin{eqnarray*}
&&\|(T-T_{h})\varphi\|_{-\frac{1}{2},\partial\Omega}=\sup\limits_{g\in H^{\frac{1}{2}}(\partial\Omega),g\neq 0}\frac{|\langle (T-T_{h})\varphi,g\rangle|}{\|g\|_{\frac{1}{2},\partial\Omega}}\\
&&~~~=\sup\limits_{g\in H^{\frac{1}{2}}(\partial\Omega),g\neq 0}\frac{|\overline{\langle g, (T-T_{h})\varphi\rangle}|}{\|g\|_{\frac{1}{2},\partial\Omega}}\\
&&~~~=\sup\limits_{g\in H^{\frac{1}{2}}(\partial\Omega),g\neq 0}\frac{|a(Ag,(A-A_{h})\varphi)|}{\|g\|_{\frac{1}{2},\partial\Omega}}\\
&&~~~=\sup\limits_{g\in H^{\frac{1}{2}}(\partial\Omega),v\in V_{h}(\Omega)}\frac{|a((A-A_{h})\varphi,Ag-v)|}{\|g\|_{\frac{1}{2},\partial\Omega}}\\
&&~~~\leq\sup\limits_{g\in H^{\frac{1}{2}}(\partial\Omega),v\in V_{h}(\Omega) }\frac{\|(A-A_{h})\varphi\|_{1,\Omega}\|Ag-v\|_{1,\Omega}}{\|g\|_{\frac{1}{2},\partial\Omega}}\\
&&~~~\lesssim h^{r}\|A\varphi\|_{1,\Omega},
\end{eqnarray*}
thus, from (\ref{s2.9}) we deduce
\begin{eqnarray*}
&&\|T-T_{h}\|_{H^{-\frac{1}{2}}(\partial\Omega)\rightarrow H^{-\frac{1}{2}}(\partial\Omega)}=\sup\limits_{\varphi\in H^{-\frac{1}{2}}(\partial\Omega),\varphi\neq 0}\frac{\|(T-T_{h})\varphi\|_{-\frac{1}{2},\partial\Omega}}{\|\varphi\|_{-\frac{1}{2},\partial\Omega}}\\
&&~~~\lesssim \sup\limits_{\varphi\in H^{-\frac{1}{2}}(\partial\Omega),\varphi\neq 0}\frac{h^{r}\|A\varphi\|_{1,\Omega}}{\|\varphi\|_{-\frac{1}{2},\partial\Omega}}\lesssim h^{r}\rightarrow 0 (h\rightarrow 0).\\
\end{eqnarray*}
Note that $T_{h}$ is a finite rank operator, thus, $T$ is a compact operator. The proof is completed. ~~~$\Box$\\

\indent In this paper, we suppose that $\{\lambda_{p}\}$ and
$\{\lambda_{p,h}\}$ are enumerations of the eigenvalues of
(\ref{s2.2}) and (\ref{s2.5}) respectively according to the same sort rule, each repeated as many
times as its multiplicity, and $\lambda=\lambda_{j}$ is the $j$th
eigenvalue with the algebraic multiplicity $q$ and the ascent
$\alpha$, $\lambda=\lambda_{j}=\lambda_{j+1}=\cdots=\lambda_{j+q-1}$.
 Since $T_{h}$ converges to $T$, $q$ eigenvalues $\lambda_{j,h}, \lambda_{j+1,h}, \cdots, \lambda_{j+q-1,h}$ of (\ref{s2.5}) will converge to $\lambda$. Let $M(\lambda)$ be the space spanned by all eigenfunctions corresponding to the
eigenvalue $\lambda$, and $M_{h}(\lambda)$ be the space spanned by all generalized eigenfunctions of (\ref{s2.5}) corresponding to the
eigenvalues $\lambda_{p,h}(p=j,j+1,\cdots,j+q-1)$. In view of the dual problem (\ref{s2.15}) and (\ref{s2.16}), the definitions of $M^*(\lambda^*)$
and $M_{h}^*(\lambda^*)$ are analogous to those of $M(\lambda)$ and $M_{h}(\lambda)$.\\
\noindent{\bf Theorem 2.1.} Let $\lambda$ and $\lambda_h$ be the $j$th eigenvalue of (\ref{s2.2}) and  (\ref{s2.5}), respectively.
Let $M(\lambda), M^*(\lambda^*)\subset H^{1+t}(\Omega)~(t\geq r)$, then
\begin{eqnarray}\label{s2.22}
|\lambda-\lambda_h|\lesssim h^{\frac{2\sigma}{\alpha}},
\end{eqnarray}
suppose that $u_{h}$ is an eigenfunction corresponding to $\lambda_{h}$,
then there exists an eigenfunction
$u$ corresponding to $\lambda$ such that
\begin{eqnarray}\label{s2.23}
&&\|u_{h}-u\|_{1,\Omega}\lesssim h^{\frac{r+\sigma}{\alpha}}+h^{\sigma},\\\label{s2.24}
&&\|u_{h}-u\|_{0,\Omega}\lesssim h^{\frac{r+\sigma}{\alpha}},\\\label{s2.25}
&&\|u_{h}-u\|_{-\frac{1}{2},\partial\Omega}\lesssim h^{\frac{r+\sigma}{\alpha}},
\end{eqnarray}
where $\sigma=\min\{m,t\}$ and the principle to determine $r$ see Lemma 2.1.\\
{\bf{Proof.}}
%
~~Since $\|T-T_{h}\|_{H^{-\frac{1}{2}}(\partial\Omega)\rightarrow  H^{-\frac{1}{2}}(\partial\Omega)}\rightarrow 0$,
from Theorem 7.3 and Theorem 7.4 in \cite{babuska} we know that
there exists an eigenfunction $u$ corresponding to $\lambda$
and
\begin{eqnarray}\label{s2.26}
&&\mid\lambda-\lambda_{h}\mid \lesssim \{\sum\limits_{l,p=j}^{j+q-1}\mid\langle(T-T_{h})\varphi_{l},\varphi_{p}^*\rangle\mid\nonumber\\
&&~~~+\|(T-T_{h}) \mid_{M(\lambda)} \|_{-\frac{1}{2},\partial\Omega}\|(T^*-T_{h}^*) \mid_{M^*(\lambda^{*})} \|_{\frac{1}{2},\partial\Omega}\}^{1/\alpha},\\\label{s2.27}
&&\|u_{h}-u\|_{-\frac{1}{2},\partial\Omega}\lesssim \|(T_{h}-T)\mid_{M(\lambda)}\|_{-\frac{1}{2},\partial\Omega}^{\frac{1}{\alpha}},
\end{eqnarray}
where $\varphi_{j},\cdots,\varphi_{j+q-1}$ are any basis for $M(\lambda)$ and $\varphi_{j}^{*},\cdots,\varphi_{j+q-1}^{*}$ are the dual basis in $M^*(\lambda^{*})$.\\
From the definitions of $A^*$ and $A_{h}$, we deduce that
\begin{eqnarray}\label{s2.28}
\mid\langle(T-T_{h})\varphi_{l},\varphi_{p}^*\rangle\mid&=&\mid a((A-A_{h})\varphi_{l}, A^*\varphi_{p}^*)\mid\nonumber\\
&=&\mid a((A-A_{h})\varphi_{l}, A^*\varphi_{p}^*-A_{h}^*\varphi_{p}^*)\mid\nonumber\\
&\lesssim& \|(A-A_{h})\varphi_{l}\|_{1,\Omega}\|A^*\varphi_{p}^*-A_{h}^*\varphi_{p}^*\|_{1,\Omega}\nonumber\\
&\lesssim& h^{2\sigma}.
\end{eqnarray}
It is easy to know that the second term on the right-hand side of  (\ref{s2.26}) is a quantity of higher order than $h^{2\sigma}$, then substituting (\ref{s2.28}) into  (\ref{s2.26}) we get (\ref{s2.22}).\\
From (\ref{s2.21}) we obtain
\begin{eqnarray}\label{s2.29}
&&\|(T-T_{h})|_{M(\lambda)}\|_{-\frac{1}{2},\partial\Omega}=\sup\limits_{f\in M(\lambda),\|f\|_{-\frac{1}{2},\partial\Omega}=1 }\|Tf-T_{h}f\|_{-\frac{1}{2},\partial\Omega}\nonumber\\
&&~~~\lesssim h^{r+\sigma} \sup\limits_{f\in M(\lambda),\|f\|_{-\frac{1}{2},\partial\Omega}=1 } \|Af\|_{1+t,\Omega}.
\end{eqnarray}
Substituting (\ref{s2.29}) into  (\ref{s2.27}) we get (\ref{s2.25}).\\
By calculation, we have
\begin{eqnarray}\label{s2.30}
u_{h}-u&=&\lambda Au-\lambda_{h} A_{h}u_{h}\nonumber\\
&=&(\lambda-\lambda_{h})Au+\lambda_{h}A(u-u_{h})+\lambda_{h}(A-A_{h})u_{h},
\end{eqnarray}
then, from (\ref{s2.22}), (\ref{s2.9}), (\ref{s2.25}) and (\ref{s2.20}) we derive
\begin{eqnarray*}
&&\|u_{h}-u\|_{1,\Omega}\leq\|(\lambda_{h}-\lambda)Au\|_{1,\Omega}+\|\lambda_{h}A(u_{h}-u)\|_{1,\Omega}+\|\lambda_{h}(A_{h}-A)u_{h}\|_{1,\Omega}\nonumber\\
&&~~~\lesssim |\lambda_{h}-\lambda|+\|A(u_{h}-u)\|_{1,\Omega}+\|(P_{h}A-A)u_{h}\|_{1,\Omega}\nonumber\\
&&~~~\lesssim h^{\frac{2\sigma}{\alpha}}+\|u_{h}-u\|_{-\frac{1}{2},\partial\Omega}+\|(P_{h}A-A)(u_{h}-u+u)\|_{1,\Omega}\nonumber\\
&&~~~\lesssim h^{\frac{2\sigma}{\alpha}}+h^{\frac{r+\sigma}{\alpha}}+\|(P_{h}A-A)(u_{h}-u)\|_{1,\Omega}+\|(P_{h}A-A)u\|_{1,\Omega}\nonumber\\
&&~~~\lesssim h^{\frac{2\sigma}{\alpha}}+h^{\frac{r+\sigma}{\alpha}}+h^{\sigma}\nonumber\\
&&~~~\lesssim h^{\frac{r+\sigma}{\alpha}}+h^{\sigma},\nonumber\\
&&\|u-u_{h}\|_{0,\Omega}\leq\|(\lambda_{h}-\lambda)Au\|_{0,\Omega}+\|\lambda_{h}A(u_{h}-u)\|_{0,\Omega}+\|\lambda_{h}(A_{h}-A)u_{h}\|_{0,\Omega}\nonumber\\
&&~~~\lesssim |\lambda_{h}-\lambda|+\|A(u_{h}-u)\|_{1,\Omega}+\|(P_{h}A-A)u_{h}\|_{0,\Omega}\nonumber\\
&&~~~\lesssim h^{\frac{2\sigma}{\alpha}}+\|u_{h}-u\|_{-\frac{1}{2},\partial\Omega}+h^{r}\|(P_{h}A-A)u_{h}\|_{1,\Omega}\nonumber\\
&&~~~\lesssim h^{\frac{2\sigma}{\alpha}}+h^{\frac{r+\sigma}{\alpha}}+h^{r+\sigma}\nonumber\\
&&~~~\lesssim h^{\frac{r+\sigma}{\alpha}}.
\end{eqnarray*}
The proof is completed.~~~$\Box$\\

\indent For the dual problem (\ref{s2.15}) and (\ref{s2.16}), we have the corresponding conclusion as Theorem 2.1.\\

\section{Local a priori error estimates}
\indent In this section, we will discuss local a priori error estimates which are a basic issue in finite element method and a basic tool for analyzing the local computational algorithm we will talk about later.\\
\indent For $D\subset G\subset \Omega$, we use the notation $D\subset\subset
G$ to mean that $dist(\partial D\backslash\partial\Omega,\partial
G\backslash\partial\Omega)>0$ (see Fig. 1).\\
\begin{center}
\includegraphics[width=0.35\textwidth]{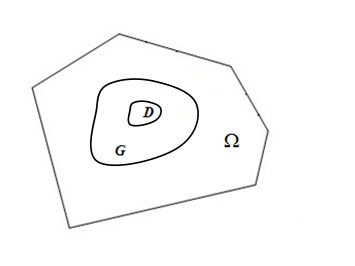}~~~
\includegraphics[width=0.35\textwidth]{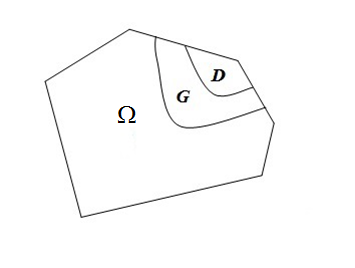}\\
\textrm{\small\bf Fig.1. Subdomains}
\end{center}
Given $G\subset\Omega$, we define $\pi_{h}(G)$ and
$V_{h}(G)$  to be the restriction of $\pi_{h}(\Omega)$ and $V_{h}(\Omega)$ to $G$,
respectively. Denote $supp~v=\overline{\{x: v(x)\not= 0\}}$, and
$$V_{h}^{0}(G)=\{v\in V_{h}(\Omega): v|_{\partial G\setminus \partial\Omega}=0\}, ~~~V_{0}^{h}(G)=\{v\in V_{h}(\Omega): (supp~v\setminus \partial\Omega)\subset\subset G\}.$$

Let $\Omega_{0}\subset \Omega$. We assume that the meshes and finite
element spaces in this paper satisfy the following
conditions (see \cite{xu1}): \\
\indent (C1)~ There exists $\gamma\geq 1$ such that
\begin{eqnarray*}
h_{\Omega}^{\gamma}\leq Ch(x), \forall x\in \Omega.
\end{eqnarray*}
\indent (C2)~ {\em Inverse Estimate.} For any $v\in
V_{h}(\Omega_{0})$,
\begin{eqnarray*}
\|v\|_{1,\Omega_{0}}\leq C(\min\limits_{x\in \Omega_{0}}
h(x))^{-1}\|v\|_{0,\Omega_{0}}.
\end{eqnarray*}
\indent (C3)~{\em Superapproximation.} For $G\subset \Omega_{0}$,
let $\omega\in C^{\infty}(\overline{\Omega})$ with
(supp $\omega\setminus\partial\Omega) \subset\subset G$. Then for any
$\varphi\in V_{h}(G)$, there exists $v\in V_{0}^{h}(G)$ such that
\begin{eqnarray*}
\|h_{G}^{-1}(\omega \varphi-v)\|_{1,G}\leq C\|\varphi\|_{1,G}.
\end{eqnarray*}
From \cite{canuto,shen} we know that (C0), (C2) and (C3) also hold for conforming spectral element.\\


\indent For some $G\subset\subset \Omega$, we consider the following mixed boundary value problem:
\begin{eqnarray}\label{s3.1}
\triangle  \psi +k^{2}n(x)  \psi&=&g~~in~ G,\nonumber\\
\psi&=&0~~ on ~\partial G\setminus \partial\Omega,\\
\frac{\partial \psi}{\partial \nu}&=&0~~ on~ \partial G \cap \partial\Omega.\nonumber
\end{eqnarray}
\indent The weak form of (\ref{s3.1}) is given by: Find
$\psi_{g} \in H^{1}_{\Gamma}(G)=\{v\in H^{1}(G): v|_{\partial G\setminus \partial\Omega}=0\}$, such that
\begin{eqnarray}\label{s3.2}
a(v,\psi_{g})=(v,g),~~~\forall v\in H^{1}_{\Gamma}(G).
\end{eqnarray}
For (\ref{s3.1}) we need the following assumption.\\
\indent {\bf R(G)}. For any $f\in L^{2}(G)$, there exists $\varphi\in H^{1+r}_{\Gamma}(G)$ satisfying 
$$a(v,\varphi)=(v,f)~~~\forall v\in H^{1}_{\Gamma}(G)$$
and
$$\|\varphi\|_{1+r,G}\leq C\|f\|_{0,G}.$$
Here, since $G$ is a local domain, we can easily control the shape of $G$ to make the intersection angle of two boundary parts is less than $\pi$, even less than or equal to $\frac{\pi}{2}$. Thus, from Theorem 2.1 and Remark 3 in \cite{babuska2}, or \cite{mitrea,taylor} we know that the above assumption {\bf R(G)} is reasonable.\\
Let
\begin{eqnarray*}
a_{0}(u,v)=\int\limits_{\Omega}\nabla u\cdot \overline{\nabla v}dx.
\end{eqnarray*}

From \cite{xu1,yang1}, after a minor modification, we have the following technical result.\\
\indent{\bf Lemma 3.1.}~~Let
$D\subset\subset\Omega_{0}\subset\Omega\subset \mathbb{R}^{2}$, and
$\omega\in C^{\infty}(\overline{\Omega})$ be a real valued function with $(supp~\omega\setminus\partial\Omega)\subset\subset\Omega_{0}$. Then
\begin{eqnarray}\label{s3.3}
a_{0}(\omega \psi,\omega \psi)\leq |a(\psi,\omega^{2}
\psi)|+C\|\psi\|_{0,\Omega_{0}}^{2}, \forall \psi\in H^{1}(\Omega).
\end{eqnarray}
{\bf{Proof.}}
By calculating, $\forall \psi\in H^{1}(\Omega)$, we have
 \begin{eqnarray*}
&&a_{0}(\omega \psi,\omega \psi)=\int_{\Omega}\sum\limits_{i=1}^{2}\frac{\partial(\omega
\psi)}{\partial x_{i}}\frac{\overline{\partial(\omega \psi)}}{\partial x_{i}}
=\int_{\Omega}\sum\limits_{i=1}^{2}(\frac{\partial \omega
}{\partial x_{i}}\psi+ \omega \frac{\partial \psi}{\partial x_{i}})\overline{(\frac{\partial \omega
}{\partial x_{i}}\psi+ \omega \frac{\partial \psi}{\partial x_{i}})}\\
&&=\int_{\Omega}\sum\limits_{i=1}^{2}(\frac{\partial \omega}{\partial x_{i}}\psi+ \omega \frac{\partial \psi}{\partial x_{i}})(\overline{\frac{\partial \omega}{\partial x_{i}}}\overline{\psi}+\overline{\omega}\overline{\frac{\partial \psi}{\partial x_{i}}})\\
&&=\int_{\Omega}\sum\limits_{i=1}^{2}(|\frac{\partial \omega}{\partial x_{i}}\psi|^{2}+2Re\{\overline{\frac{\partial \omega}{\partial x_{i}}}\frac{\partial \psi}{\partial x_{i}}\omega \overline{\psi}\}+|\omega \frac{\partial \psi}{\partial x_{i}}|^{2}),
\end{eqnarray*}
and
\begin{eqnarray*}
&&a_{0}(\psi,\omega^{2} \psi)=\int_{\Omega}\sum\limits_{i=1}^{2}\frac{\partial
\psi}{\partial x_{i}}\overline{\frac{\partial(\omega^2 \psi)}{\partial x_{i}}} =\int_{\Omega}\sum\limits_{i=1}^{2}\frac{\partial
\psi}{\partial x_{i}}\overline{(2\omega \frac{\partial\omega}{\partial x_{i}}\psi
+\omega^{2} \frac{\partial \psi}{\partial x_{i}})}\\
&&=\int_{\Omega}\sum\limits_{i=1}^{2}\frac{\partial \psi}{\partial x_{i}}(2\overline{\omega}\overline{\frac{\partial\omega}{\partial x_{i}}}\overline{\psi}+\overline{\omega^{2}}\overline{\frac{\partial \psi}{\partial x_{i}}})=\int_{\Omega}\sum\limits_{i=1}^{2}(2\frac{\partial \psi}{\partial x_{i}}\overline{\frac{\partial\omega}{\partial x_{i}}}\overline{\omega}\overline{\psi}+|\omega \frac{\partial \psi}{\partial x_{i}}|^{2}).
\end{eqnarray*}
Noting that $\omega$ is a real valued function and comparing the above two relations, we have
\begin{eqnarray}\label{s3.4}
a_{0}(\omega \psi,\omega \psi)&=&Re\{a_{0}(\psi,\omega^{2} \psi)\}+
\int_{\Omega}\sum\limits_{i=1}^{2}|\frac{\partial \omega}{\partial x_{i}} \psi|^{2}\nonumber\\
&=&Re\{a(\psi,\omega^{2} \psi)+k^{2}(n\psi,\omega^{2} \psi)\}+\int_{\Omega}\sum\limits_{i=1}^{2}|\frac{\partial \omega}{\partial x_{i}} \psi|^{2}\nonumber\\
&\leq&|a(\psi,\omega^{2} \psi)|+K(n_{1}\psi,\omega^{2} \psi)+\int_{\Omega}\sum\limits_{i=1}^{2}|\frac{\partial
\omega}{\partial x_{i}} \psi|^{2},
\end{eqnarray}
where $K$ is a positive constant that is large enough.\\
Since $\omega\in C^{\infty}(\overline{\Omega})$ , $(supp~\omega\setminus\partial\Omega)\subset\subset\Omega_{0}$,
\begin{eqnarray*}
\int_{\Omega}\sum\limits_{i=1}^{2}|\frac{\partial
\omega}{\partial x_{i}} \psi|^{2}\leq C\|\psi\|_{0,\Omega_{0}}^{2},
\end{eqnarray*}
which together with (\ref{s3.4}) yields (\ref{s3.3}).
The proof is completed.~~~$\Box$\\

Using the proof method in \cite{xu1,yang1} we prove the following conclusions.\\
\indent{\bf Lemma 3.2.}~~Suppose that $D\subset\subset\Omega_{0}$ and (C1), (C2)
and (C3) hold. If $f\in H^{-1}(\Omega)$ and $\psi\in V_{h}(\Omega_{0})$ satisfies
\begin{eqnarray}\label{s3.5}
a(\psi,v)=f(v),~~~\forall v \in V_{0}^{h}(\Omega_{0}),
\end{eqnarray}
then
\begin{eqnarray}\label{s3.6}
\|\psi\|_{1,D}\lesssim \|\psi\|_{0,\Omega_{0}}+\|f\|_{-1,\Omega}.
\end{eqnarray}
 \indent{\bf Proof.}~~Let $p\geq 2\gamma-1$ be an integer, and let
 $$D\subset\subset \Omega_{p}\subset\subset \Omega_{p-1}\subset\subset\cdots\subset\subset \Omega_{1}\subset\subset \Omega_{0}.$$
 Choose $D_{1}\subset \Omega$ satisfying $D\subset\subset D_{1}\subset\subset \Omega_{p}$
and $\omega\in C^{\infty}(\overline{\Omega})$ such that
$\omega\equiv 1$ on $\overline{D_{1}}$ and $(supp~\omega\setminus\partial\Omega)\subset\subset\Omega_{p}$. Then, from
(C3), there exists $v\in V_{0}^{h}(\Omega_{p})$ such that
\begin{eqnarray*}
\|\omega^{2} \psi-v\|_{1,\Omega_{p}}\lesssim
h_{\Omega_{0}}\|\psi\|_{1,\Omega_{p}},
\end{eqnarray*}
so we have
\begin{eqnarray}\label{s3.7}
|a(\psi, \omega^{2} \psi-v)|\lesssim h_{\Omega_{0}}\|\psi\|_{1,\Omega_{p}}^{2},
\end{eqnarray}
and from the trace theorem we get
\begin{eqnarray}\label{s3.8}
&&|f(v)|\lesssim
\|f\|_{-1,\Omega}\|v\|_{1,\Omega}\lesssim
\|f\|_{-1,\Omega}\|v\|_{1,\Omega_{p}}\nonumber\\
&&\lesssim
\|f\|_{-1,\Omega}(h_{\Omega_{0}}\|\psi\|_{1,\Omega_{p}}+\|\omega
\psi\|_{1,\Omega}).
\end{eqnarray}
Since $v\in V_{0}^{h}(\Omega_{p})\subset  V_{0}^{h}(\Omega_{0})$,
(\ref{s3.5}) implies
\begin{eqnarray}\label{s3.9}
a(\psi, \omega^{2} \psi)=a(\psi, \omega^{2} \psi-v)+f(v).
\end{eqnarray}
It follows from (\ref{s3.3}), (\ref{s3.9}), (\ref{s3.7}) and (\ref{s3.8}) that
\begin{eqnarray*}
&&\|\omega \psi\|_{1,\Omega}^{2}\lesssim a_{0}(\omega \psi,\omega \psi)\lesssim |a(\psi,\omega^{2} \psi)|+\|\psi\|_{0,\Omega_{0}}^{2}\\
&&\leq |a(\psi,\omega^{2} \psi-v)|+\|\psi\|_{0,\Omega_{0}}^{2}+|f(v)|\\
&&\lesssim
h_{\Omega_{0}}\|\psi\|_{1,\Omega_{p}}^{2}+\|\psi\|_{0,\Omega_{0}}^{2}+\|f\|_{-1,\Omega}(h_{\Omega_{0}}\|\psi\|_{1,\Omega_{p}}+\|\omega
\psi\|_{1,\Omega}),
\end{eqnarray*}
thus
\begin{eqnarray}\label{s3.10}
 \|\psi\|_{1,D}\lesssim
 h_{\Omega_{0}}^{1/2}\|\psi\|_{1,\Omega_{p}}+\|\psi\|_{0,\Omega_{0}}+\|f\|_{-1,\Omega}.
\end{eqnarray}
Similarly, we can get
\begin{eqnarray}\label{s3.11}
 \| \psi\|_{1,\Omega_{j}}\lesssim
 h_{\Omega_{0}}^{1/2}\|\psi\|_{1,\Omega_{j-1}}+\|\psi\|_{0,\Omega_{0}}+\|f\|_{-1,\Omega},
 ~~~j=1,2,\cdots,p.
\end{eqnarray}
By using (\ref{s3.10}) and (\ref{s3.11}), we get from (C1) and (C2) that
\begin{eqnarray*}
&&\| \psi\|_{1,D}\lesssim
 h_{\Omega_{0}}^{(p+1)/2}\|\psi\|_{1,\Omega_{0}}+\|\psi\|_{0,\Omega_{0}}+\|f\|_{-1,\Omega}\\
 &&\lesssim
 h_{\Omega_{0}}^{(p+1)/2}\|(\min\limits_{x\in \Omega_{0}}
h(x))^{-1}\psi\|_{0,\Omega_{0}}+\|\psi\|_{0,\Omega_{0}}+\|f\|_{-1,\Omega}\\
 &&\lesssim  \|\psi\|_{0,\Omega_{0}}+\|f\|_{-1,\Omega}.
\end{eqnarray*}
This completes the proof.~~~$\Box$\\

\indent{\bf Theorem 3.1.} Suppose that $\psi\in H^{1}(\Omega)$,
$D\subset\subset\Omega_{0}$, (C0), (C1), (C2), (C3) and  {\bf R($\Omega_{0}$)}  hold. Then
\begin{eqnarray}\label{s3.14}
\|P_{h}\psi\|_{1,D}\lesssim
\|\psi\|_{1,\Omega_{0}}+\|P_{h}\psi\|_{0,\Omega_{0}}.
\end{eqnarray}
 \indent{\bf Proof.}~~Let $P_{h}^{\Omega_{0}}:H^{1}(\Omega_{0})\to
 V_{h}(\Omega_{0})$ be the projection defined by
\begin{eqnarray}\label{s3.14a}
a(\psi-P_{h}^{\Omega_{0}}\psi, v)=0,~~~\forall v\in V_{h}(\Omega_{0}).
\end{eqnarray}
Similar to (\ref{s2.37}) we can prove that $\|P_{h}^{\Omega_{0}}\psi\|_{1,\Omega_{0}}\lesssim \|\psi\|_{1,\Omega_{0}}$ for any $\psi\in H^{1}(\Omega_{0})$.\\
\indent Choose $D_{1}\subset\subset\Omega$ such that  $D\subset\subset
D_{1}\subset\subset\Omega_{0}$, $\omega\in
C^{\infty}(\overline{\Omega})$, $\omega\equiv 1$ on
$\overline{D}_{1}$, $(supp~\omega\setminus\partial\Omega)\subset\subset\Omega_{0}$. Let
$\widetilde{\psi}=\omega \psi$, then $\forall v\in V_{0}^{h}(D_{1})$,
\begin{eqnarray*}
&&a(P_{h}^{\Omega_{0}}\widetilde{\psi}-P_{h}\psi,
v)=a(P_{h}^{\Omega_{0}}(\omega \psi)-\omega \psi+\omega \psi-P_{h}\psi, v)\\
&&~~~=a(P_{h}^{\Omega_{0}}(\omega \psi)-\omega \psi, v)+a(\omega \psi-P_{h}\psi,
v)=0.
\end{eqnarray*}
Thus, from Lemma 3.2, we have
\begin{eqnarray*}
\|P_{h}^{\Omega_{0}}\widetilde{\psi}-P_{h}\psi\|_{1,D}\lesssim
\|P_{h}^{\Omega_{0}}\widetilde{\psi}-P_{h}\psi\|_{0,D_{1}},
\end{eqnarray*}
then, we derive
\begin{eqnarray*}
&&\|P_{h}\psi\|_{1,D}\leq \|P_{h}^{\Omega_{0}}\widetilde{\psi}\|_{1,D}+
\|P_{h}^{\Omega_{0}}\widetilde{\psi}-P_{h}\psi\|_{1,D}\\
&&~~~\lesssim \|P_{h}^{\Omega_{0}}\widetilde{\psi}\|_{1,D}+
\|P_{h}^{\Omega_{0}}\widetilde{\psi}-P_{h}\psi\|_{0,D_{1}}\\
&&~~~\lesssim \|P_{h}^{\Omega_{0}}\widetilde{\psi}\|_{1,D_{1}}+
\|P_{h}\psi\|_{0,D_{1}}\\
&&~~~\lesssim \|\widetilde{\psi}\|_{1,\Omega_{0}}+
\|P_{h}\psi\|_{0,D_{1}}\\
&&~~~\lesssim \|\psi\|_{1,\Omega_{0}}+ \|P_{h}\psi\|_{0,\Omega_{0}}.
\end{eqnarray*}
The proof is completed.~~~$\Box$\\

Using Theorem 3.1 and the proof method of Theorems 3.4-3.5 in \cite{yang1} we can prove the following local estimates.\\
\indent{\bf Theorem 3.2.} Suppose that $\psi\in H^{1}(\Omega)$,
$D\subset\subset\Omega_{0}$, (C0), (C1), (C2) and (C3) hold. Then
\begin{eqnarray}\label{s3.15}
\|\psi-P_{h}\psi\|_{1,D}\lesssim \inf\limits_{v\in
V_{h}(\Omega)}\|\psi-v\|_{1,\Omega_{0}}+\|\psi-P_{h}\psi\|_{0,\Omega}.
\end{eqnarray}

\indent{\bf Theorem 3.3.} Under the assumptions of Theorem 3.2, let
$(\lambda_{h}, u_{h})$ be the $j$th eigenpair of
 (\ref{s2.5}), and $\lambda$ be the $j$th
eigenvalue of (\ref{s2.2}).
Then there exists an eigenfunction $u$ corresponding to $\lambda$ such that the following error
estimate holds:
\begin{eqnarray}\label{s3.16}
\|u-u_{h}\|_{1,D}&\lesssim &\inf\limits_{v\in
V_{h}(\Omega)}\|u-v\|_{1,\Omega_{0}}+\|u-P_{h}u\|_{0,\Omega}\nonumber\\
&~~~&+\|\lambda u-\lambda_{h}u_{h}\|_{-\frac{1}{2},\partial\Omega}.
\end{eqnarray}

\section{Two-grid discretizations for the Stekloff eigenvalue problem}
In this section, we present two kinds of two-grid discretizations for the Stekloff eigenvalue problem (\ref{s2.1}).\\
\noindent{\bf Scheme 1.}(Two-grid scheme I)\\
\noindent{\bf Step 1.} Solve (\ref{s2.5}) on a coarse grid $\pi_{H}$: Find
$(\lambda_{H}, u_{H})\in \mathbb{C}\times V_{H}(\Omega)$, $u_{H}\neq 0$ such that
\begin{eqnarray*}
a(u_{H},v)=-\lambda_{H} \langle u_{H},v\rangle,~~~\forall v\in V_{H}(\Omega).
\end{eqnarray*}
Let $\lambda_{H}^*=\overline{\lambda_{H}}$, and find $u_{H}^{*}\in M_{H}^*(\lambda^*)$ such that
$|\langle u_{H},u_{H}^{*}\rangle|$ has a positive lower bound uniformly with respect to $H$.\\
\noindent{\bf Step 2.} Solve two linear boundary value problems on a fine grid  $\pi_{w}(w<H)$:
Find $u^{w}\in V_{w}(\Omega)$ such that
\begin{eqnarray*}
a(u^{w}, v)=-\lambda_{H}\langle u_{H},v\rangle,~~~\forall v\in V_{w}(\Omega);
\end{eqnarray*}
find $u^{w*}\in V_{w}(\Omega)$ such that
\begin{eqnarray*}
a(v, u^{w*})=-\lambda_{H}\langle v, u_{H}^{*}\rangle,~~~\forall v\in V_{w}(\Omega).
\end{eqnarray*}
\noindent{\bf Step 3.} Compute the generalized Rayleigh quotient
$$\lambda^{w}=-\frac{a(u^{w},u^{w*})}{\langle u^{w},u^{w*}\rangle}.$$

\noindent{\bf Scheme 2.}(Two-grid scheme II)\\
\noindent{\bf Step 1.} The same as Step 1 of Scheme 1.\\
\noindent{\bf Step 2.} Solve two linear boundary value problems on a fine grid  $\pi_{w}$:
Find $u^{w}\in V_{w}(\Omega)$ such that
\begin{eqnarray*}
(\nabla u^{w}, \nabla v)+(u^{w}, v)=-\lambda_{H}\langle u_{H},v\rangle+((k^{2}n+1)u_{H}, v),~~~\forall v\in V_{w}(\Omega);
\end{eqnarray*}
find $u^{w*}\in V_{w}(\Omega)$ such that
\begin{eqnarray*}
(\nabla v, \nabla u^{w*})+(v, u^{w*})=-\lambda_{H}\langle v, u_{H}^{*}\rangle+(v, (k^{2}n+1)u_{H}^{*}),~~~\forall v\in V_{w}(\Omega).
\end{eqnarray*}
\noindent{\bf Step 3.} Compute the generalized Rayleigh quotient
$$\lambda^{w}=-\frac{a(u^{w},u^{w*})}{\langle u^{w},u^{w*}\rangle}.$$

\indent\noindent{\bf Remark.} Let $u_{H}^{-}$ be the orthogonal projection of $u_{H}$ to $M_{H}^*(\lambda^*)$ in the sense of the inner product $\langle \cdot, \cdot\rangle$, and $u_{H}^*=\frac{u_{H}^{-}}{\|u_{H}^{-}\|_{0,\partial\Omega}}$, then when $H$ is small enough $|\langle u_{H}, u_{H}^*\rangle|$
has a positive lower bound uniformly with respect to $H$. One can refer to \cite{yang2016} for the proof of this conclusion. Therefore, $u_{H}^{*}$ in Step 1 of Schemes 1 and 2 can be obtained in this way.\\


\indent{\bf Lemma 4.1.} Let $(\lambda, u)$ and $(\lambda^*, u^*)$ be the eigenpair of (\ref{s2.2}) and  (\ref{s2.15}), respectively.
Then, for any $v, v^*\in H^{1}(\Omega)$, $\langle v, v^*\rangle\neq 0$, the generalized Rayleigh quotient satisfies
\begin{eqnarray}\label{s4.1}
-\frac{a(v, v^*)}{\langle v, v^*\rangle}-\lambda=-\frac{a(v-u, v^*-u^*)}{\langle v, v^*\rangle}-\lambda\frac{\langle v-u, v^*-u^*\rangle}{\langle v, v^*\rangle}.
\end{eqnarray}
\indent {\bf{Proof.}} From (\ref{s2.2}) and (\ref{s2.15}), and by a simple calculation, we have
\begin{eqnarray*}
-a(v-u, v^*-u^*)-\lambda\langle v-u, v^*-u^*\rangle=-a(v, v^*)-\lambda\langle v, v^*\rangle,
\end{eqnarray*}
dividing both sides by $\langle v, v^*\rangle$ we obtain the desired result.~~~$\Box$\\

\indent{\bf Theorem 4.1.} Let $(\lambda_{H}, u_{H})$, $(\lambda_{H}^*, u_{H}^*)$, $u^{w}$,  $u^{w*}$, and $\lambda^{w}$ be obtained
by Scheme 1. Let $M(\lambda), M^*(\lambda^*)\subset H^{1+t}(\Omega)(t\geq r)$, then there exists an eigenfunction $u\in M(\lambda)$ and an eigenfunction $u^*\in M^*(\lambda^*)$ such that
\begin{eqnarray}\label{s4.2}
&&\|u^w-u\|_{1,\Omega}\lesssim H^{\frac{r+\sigma}{\alpha}}+w^{\sigma},\\\label{s4.3}
&&\|u^{w*}-u^*\|_{1,\Omega}\lesssim H^{\frac{r+\sigma}{\alpha}}+w^{\sigma},\\\label{s4.4}
&&\|u^w-u\|_{0,\Omega}\lesssim H^{\frac{r+\sigma}{\alpha}},\\\label{s4.3a}
&&\|u^{w*}-u^*\|_{0,\Omega}\lesssim H^{\frac{r+\sigma}{\alpha}};
\end{eqnarray}
further, assume that the ascent of $\lambda$ is equal to 1, then
\begin{eqnarray}\label{s4.4a}
&&|\lambda^{w}-\lambda|\lesssim H^{2r+2\sigma}+w^{2\sigma},
\end{eqnarray}
where $\sigma=\min\{m,t\}$.\\
\indent {\bf{Proof.}} Let $u\in M(\lambda)$ such that $u_{H}-u$ and $\lambda_{H}-\lambda$ satisfy Theorem 2.1.
From (\ref{s2.12}) we get $u=-\lambda Au$, and from the definition of $A_{w}$ and Step 2 of Scheme 1 we get $u^{w}=-\lambda_{H} A_{w}u_{H}$.
Then,
\begin{eqnarray*}
&&\|\lambda_{H} A_{w}u_{H}-\lambda A_{w}u\|_{1,\Omega}^2\lesssim  \widetilde{a}(\lambda_{H} A_{w}u_{H}-\lambda A_{w}u, \lambda_{H} A_{w}u_{H}-\lambda A_{w}u)\\
&&~~~= a(\lambda_{H} A_{w}u_{H}-\lambda A_{w}u, \lambda_{H} A_{w}u_{H}-\lambda A_{w}u)+K(\lambda_{H} A_{w}u_{H}-\lambda A_{w}u, \lambda_{H} A_{w}u_{H}-\lambda A_{w}u)\\
&&~~~= \langle \lambda_{H}u_{H}-\lambda u, \lambda_{H} T_{w}u_{H}-\lambda T_{w}u \rangle+K(\lambda_{H} A_{w}u_{H}-\lambda A_{w}u, \lambda_{H} A_{w}u_{H}-\lambda A_{w}u)\\
&&~~~\leq \|\lambda_{H}u_{H}-\lambda u\|_{-\frac{1}{2},\partial\Omega}\|\lambda_{H} T_{w}u_{H}-\lambda T_{w}u \|_{\frac{1}{2},\partial\Omega}+K\|\lambda_{H} A_{w}u_{H}-\lambda A_{w}u \|_{0,\Omega}^2\\
&&~~~\leq (\|\lambda_{H}u_{H}-\lambda u\|_{-\frac{1}{2},\partial\Omega}+K\|A_{w}(\lambda_{H} u_{H}-\lambda u) \|_{1,\Omega})\|\lambda_{H} A_{w}u_{H}-\lambda A_{w}u \|_{1,\Omega}\\
&&~~~\lesssim (\|\lambda_{H}u_{H}-\lambda u\|_{-\frac{1}{2},\partial\Omega}+K\|\lambda_{H} u_{H}-\lambda u\|_{-\frac{1}{2},\partial\Omega})\|\lambda_{H} A_{w}u_{H}-\lambda A_{w}u \|_{1,\Omega}\\
&&~~~\lesssim \|\lambda_{H}u_{H}-\lambda u\|_{-\frac{1}{2},\partial\Omega}\|\lambda_{H} A_{w}u_{H}-\lambda A_{w}u \|_{1,\Omega},
\end{eqnarray*}
so we have
\begin{eqnarray*}
\|\lambda_{H} A_{w}u_{H}-\lambda A_{w}u\|_{1,\Omega}\lesssim \|\lambda_{H}u_{H}-\lambda u\|_{-\frac{1}{2},\partial\Omega}.
\end{eqnarray*}
Therefore, from Lemma 2.2, Theorem 2.1 and Lemma 2.4 we get
\begin{eqnarray*}
&&\|u^{w}-u\|_{1,\Omega}= \|-\lambda_{H}A_{w}u_{H}+\lambda Au\|_{1,\Omega}\\
&&~~~~~~\leq \|\lambda_{H}A_{w}u_{H}-\lambda A_{w}u\|_{1,\Omega}+\|\lambda A_{w}u-\lambda Au\|_{1,\Omega}\\
&&~~~~~~\lesssim\|\lambda_{H}u_{H}-\lambda u\|_{-\frac{1}{2},\partial\Omega}+\|A_{w}(\lambda u)-A(\lambda u)\|_{1,\Omega}\\
&&~~~~~~\lesssim H^{\frac{2\sigma}{\alpha}}+H^{\frac{r+\sigma}{\alpha}}+w^{\sigma}\lesssim H^{\frac{r+\sigma}{\alpha}}+w^{\sigma},
\end{eqnarray*}
\begin{eqnarray*}
&&\|u^{w}-u\|_{0,\Omega}= \|-\lambda_{H}A_{w}u_{H}+\lambda Au\|_{0,\Omega}\\
&&~~~~~~\leq \|\lambda_{H}A_{w}u_{H}-\lambda A_{w}u\|_{0,\Omega}+\|\lambda A_{w}u-\lambda Au\|_{0,\Omega}\\
&&~~~~~~\leq \|\lambda_{H}A_{w}u_{H}-\lambda A_{w}u\|_{1,\Omega}+\|\lambda A_{w}u-\lambda Au\|_{0,\Omega}\\
&&~~~~~~\lesssim\|\lambda_{H}u_{H}-\lambda u\|_{-\frac{1}{2},\partial\Omega}+\|A_{w}(\lambda u)-A(\lambda u)\|_{0,\Omega}\\
&&~~~~~~\lesssim H^{\frac{2\sigma}{\alpha}}+H^{\frac{r+\sigma}{\alpha}}+w^{r+\sigma}\lesssim H^{\frac{r+\sigma}{\alpha}}.
\end{eqnarray*}
Similarly we can prove (\ref{s4.3}) and (\ref{s4.3a}).\\
\indent From Lemma 4.1 we have
\begin{eqnarray}\label{s4.5}
\lambda^{w}-\lambda=-\frac{a(u^{w}-u, u^{w*}-u^*)}{\langle u^{w}, u^{w*}\rangle}-\lambda\frac{\langle u^{w}-u, u^{w*}-u^*\rangle}{\langle u^{w}, u^{w*}\rangle}.
\end{eqnarray}
Note that $u_{H}$ and $u^{w}$ just approximate the same eigenfunction $u$, $u_{H}^*$ and $u^{w*}$ approximate the same eigenfunction $u^*$,
and $\langle u_{H}, u_{H}^*\rangle$ has a positive lower bound uniformly with respect to $H$.
Hence, from
 \begin{eqnarray*}
\langle u^{w}, u^{w*}\rangle=(\langle u^{w}, u^{w*}\rangle-\langle u, u^*\rangle)+(\langle u, u^*\rangle-\langle u_{H}, u_{H}^*\rangle)+\langle u_{H}, u_{H}^*\rangle,
\end{eqnarray*}
we know that $|\langle u^{w}, u^{w*}\rangle|$ has a positive lower bound uniformly. Therefore, from (\ref{s4.5}) we get
 \begin{eqnarray}\label{s4.6}
|\lambda^{w}-\lambda|\lesssim\|u^{w}-u\|_{1,\Omega}\|u^{w*}-u^*\|_{1,\Omega}.
\end{eqnarray}
Substituting (\ref{s4.2}) and (\ref{s4.3}) with $\alpha=1$ into (\ref{s4.6}) we get (\ref{s4.4a}). The proof is completed.~~~$\Box$\\

\indent{\bf Theorem 4.2.} Let $(\lambda_{H}, u_{H})$, $(\lambda_{H}^*, u_{H}^*)$, $u^{w}$,  $u^{w*}$, and $\lambda^{w}$ be obtained
by Scheme 2. Let $M(\lambda), M^*(\lambda^*)\subset H^{1+t}(\Omega)(t\geq r)$, then there exists an eigenfunction $u\in M(\lambda)$ and an eigenfunction $u^*\in M^*(\lambda^*)$ such that
\begin{eqnarray}\label{s4.7}
&&\|u^w-u\|_{1,\Omega}\lesssim H^{\frac{r+\sigma}{\alpha}}+w^\sigma,\\\label{s4.8}
&&\|u^{w*}-u^*\|_{1,\Omega}\lesssim H^{\frac{r+\sigma}{\alpha}}+w^\sigma;
\end{eqnarray}
further, assume that the ascent of $\lambda$ is equal to 1, then
\begin{eqnarray}\label{s4.9}
|\lambda^{w}-\lambda|\lesssim H^{2r+2\sigma}+w^{2\sigma},
\end{eqnarray}
where $\sigma=\min\{m,t\}$.\\
\indent {\bf{Proof.}}~
We write (\ref{s2.2}) as follows:
\begin{eqnarray}\label{s4.12}
(\nabla u, \nabla v)+(u, v)=-\lambda\langle u,v\rangle+((k^{2}n+1)u, v),~~~\forall v\in H^1(\Omega).
\end{eqnarray}
We regard $\lambda$ and $u$ on the right-hand side of (\ref{s4.12}) as fixed, and due to the ellipticity of the left-hand side of (\ref{s4.12}) we can define Ritz projection $\widehat{P}_{w}u$ of $u$ onto $V_{w}(\Omega)$, namely,
\begin{eqnarray}\label{s4.121}
(\nabla \widehat{P}_{w}u, \nabla v)+(\widehat{P}_{w}u, v)=-\lambda\langle u,v\rangle+((k^{2}n+1)u, v),~~~\forall v\in V_{w}(\Omega),
\end{eqnarray}
then, the first equation in Step 2 of Scheme 2 minus (\ref{s4.121}) and taking $v=u^w-\widehat{P}_{w}u$ in the resulting equation we derive
 \begin{eqnarray}\nonumber
&&(\nabla (u^w-\widehat{P}_{w}u), \nabla (u^w-\widehat{P}_{w}u))+(u^w-\widehat{P}_{w}u, u^w-\widehat{P}_{w}u)\\\label{s4.122}
&&~~~~=-\langle\lambda_{H} u_{H}-\lambda u, u^w-\widehat{P}_{w}u\rangle+((k^{2}n+1)(u_{H}-u), u^w-\widehat{P}_{w}u),
\end{eqnarray}
thus we get
 \begin{eqnarray*}
\|u^w-\widehat{P}_{w}u\|_{1,\Omega}\lesssim |\lambda_{H}-\lambda|+\|u_{H}-u\|_{-\frac{1}{2}, \partial\Omega}+\|u_{H}-u\|_{0, \Omega}.
\end{eqnarray*}
Hence, from the triangle inequality, Theorem 2.1 and the error estimate of Ritz projection we deduce that
\begin{eqnarray*}
\|u^w-u\|_{1,\Omega}\leq \|u^w-\widehat{P}_{w}u\|_{1,\Omega}+\|u-\widehat{P}_{w}u\|_{1,\Omega}\lesssim H^{\frac{r+\sigma}{\alpha}}+ w^\sigma.
\end{eqnarray*}
Similarly, we can prove  (\ref{s4.8}). And from the proof of  (\ref{s4.4a}) we can obtain  (\ref{s4.9}).~~~$\Box$

\section{A local finite element scheme and the its error estimate}
\indent 
In this section, base on the two-grid discretizations and referring to Algorithm B0 in \cite{dai} we establish a local computational scheme as follows.\\
\indent Let $\pi_{H}(\Omega)$ be a shape-regular grid of size $H\in (0,1)$, $D\subset \Omega$ be a subdomain
which contains a singular point, and $\Omega_{0}$ be a slightly larger subdomain containing $D$ (namely $D\subset\subset \Omega_{0}$).
Let $\pi_{w}(\Omega)$ be a refined mesoscopic shape-regular grid (from $\pi_{H}(\Omega)$)
and $\pi_{h}(\Omega_{0})$ a locally refined grid (from $\pi_{w}(\Omega)$)) that satisfy
$h\ll w\ll H$. In our discussion, we shall use an auxiliary fine grid $\pi_{h}(\Omega)$ which is globally defined.\\
\noindent{\bf Scheme 3.}(A local scheme)\\
\noindent{\bf Step 1.} The same as Step 1 of Scheme 1.\\
\noindent{\bf Step 2.} Solve two linear boundary value problems on a globally mesoscopic grid  $\pi_{w}(\Omega)$:
Find $u^{w}\in V_{w}(\Omega)$ such that
\begin{eqnarray*}
a(u^{w}, v)=-\lambda_{H}\langle u_{H},v\rangle,~~~\forall v\in V_{w}(\Omega);
\end{eqnarray*}
find $u^{w*}\in V_{w}(\Omega)$ such that
\begin{eqnarray*}
a(v, u^{w*})=-\lambda_{H}\langle v, u_{H}^{*}\rangle,~~~\forall v\in V_{w}(\Omega).
\end{eqnarray*}
\noindent{\bf Step 3.} Solve two linear boundary value problems on a locally fine grid $\pi_{h}(\Omega_{0})$: Find $e^{h}\in V_{h}^{0}(\Omega_{0})$
such that
\begin{eqnarray*}
a(e^{h}, v)=-\lambda_{H}\langle u_{H},v\rangle-a(u^{w}, v),~~~\forall v\in V_{h}^{0}(\Omega_{0});
\end{eqnarray*}
find $e^{h*}\in V_{h}^{0}(\Omega_{0})$
such that
\begin{eqnarray*}
a(v,e^{h*})=-\lambda_{H}\langle v, u_{H}^*\rangle-a(v, u^{w*}),~~~\forall v\in V_{h}^{0}(\Omega_{0});
\end{eqnarray*}
\noindent{\bf Step 4.} Set
\[u^{w,h}=
\left \{
\begin{array}{ll}
u^{w}+e^{h}&~~~on~ \overline{\Omega}_{0},\\
u^{w}&~~~in~ \Omega\setminus\overline{\Omega}_{0},
\end{array}
\right.
\]
\[u^{w,h*}=
\left \{
\begin{array}{ll}
u^{w*}+e^{h*}&~~~on~ \overline{\Omega}_{0},\\
u^{w*}&~~~in~ \Omega\setminus\overline{\Omega}_{0},
\end{array}
\right.
\]
and compute the generalized Rayleigh quotient
$$\lambda^{w,h}=-\frac{a(u^{w,h},u^{w,h*})}{ \langle u^{w,h},u^{w,h*}\rangle},~~\lambda^{w,h*}=\overline{\lambda^{w,h}}.$$

\indent Next, we shall analyze the error estimation of Scheme 3.\\
\indent {\bf Theorem 5.1} Assume that $u^{w,h}$, $u^{w,h*}$ and $\lambda^{w,h}$ are obtained by Scheme 3 and the assumption $R(\Omega_{0})$ holds.
If $M(\lambda), M^*(\lambda^*)\subset H^{1}(\Omega)\cap H^{1+r}(\Omega)\cap H^{2}(\Omega\setminus\overline{D})$,
then there exists an eigenfunction $u$ corresponding to $\lambda$ and an eigenfunction $u^*$ corresponding to $\lambda^*$ such that
\begin{eqnarray}\label{s5.1}
\|u-u^{w,h}\|_{1,\Omega}\lesssim h^{r}+w+H^{\frac{2r}{\alpha}},\\\label{s5.1a}
\|u^*-u^{w,h*}\|_{1,\Omega}\lesssim h^{r}+w+H^{\frac{2r}{\alpha}},
\end{eqnarray}
further, assume that the ascent of $\lambda$ is equal to 1, then
\begin{eqnarray}\label{s5.2}
|\lambda-\lambda^{w,h}|\lesssim h^{2r}+w^2+H^{\frac{4r}{\alpha}},
\end{eqnarray}
where the principle to determine $r$ see Lemma 2.1.\\
\indent {\bf Proof.} Let $u\in M(\lambda)$ such that $u_{H}-u$ and $\lambda_{H}-\lambda$ satisfy Theorem 2.1.
Because of
$$\|u-u^{w,h}\|_{1,\Omega}\leq \|u-P_{h}u\|_{1,\Omega}+\|P_{h}u-u^{w,h}\|_{1,\Omega}$$
and Lemma 2.2, we only need to estimate $\|P_{h}u-u^{w,h}\|_{1,\Omega}$.\\
Choose $G\subset \Omega$ satisfying $D\subset\subset G\subset\subset \Omega_{0}$. Since
$$\|P_{h}u-u^{w,h}\|_{1,\Omega}\leq\|P_{h}u-u^{w,h}\|_{1,D}+\|P_{h}u-u^{w,h}\|_{1,G\setminus \overline{D}}+\|P_{h}u-u^{w,h}\|_{1,\Omega\setminus  \overline{G}},$$
we will estimate $\|P_{h}u-u^{w,h}\|_{1,D}$, $\|P_{h}u-u^{w,h}\|_{1,G\setminus \overline{D}}$,
and $\|P_{h}u-u^{w,h}\|_{1,\Omega\setminus\overline{G}}$  one by one. For these purpose, we take $F\subset \Omega$ such that $D\subset\subset F\subset\subset G\subset\subset \Omega_{0}$.\\
\indent First, from the equation
\begin{eqnarray*}
a(P_{w}u-u^{w},v)=-\lambda \langle u,v\rangle+\lambda_{H}\langle u_{H},v\rangle, ~~~\forall v\in V_{w}(\Omega),
\end{eqnarray*}
and
\begin{eqnarray}\label{s5.3}
-\lambda \langle u,v\rangle+\lambda_{H}\langle u_{H},v\rangle=(\lambda_{H}-\lambda)\langle u,v\rangle+\lambda_{H}\langle u_{H}-u,v\rangle,
\end{eqnarray}
we have by using the $H^1$-coerciveness of $\tilde{a}(\cdot,\cdot)$ that
\begin{eqnarray}\label{s5.4}
\|P_{w}u-u^{w}\|_{1,\Omega}\leq |\lambda-\lambda_{H}|+\|u-u_{H}\|_{-\frac{1}{2},\partial\Omega}+K\|P_{w}u-u^{w}\|_{0,\Omega}.
\end{eqnarray}
\indent From Step 3 of Scheme 3 we have the following identity
\begin{eqnarray}\label{s5.5}
a(u^{w,h}-P_{h}u,v)=-\lambda_{H}\langle u_{H},v\rangle+\lambda \langle u,v\rangle, ~~~\forall v\in V_{h}^{0}(\Omega_{0}),
\end{eqnarray}
then, from (\ref{s5.3}) and Lemma 3.2 we can derive
\begin{eqnarray}\label{s5.6}
\|P_{h}u-u^{w,h}\|_{1,D}\lesssim \|P_{h}u-u^{w,h}\|_{0,\Omega_{0}}+|\lambda-\lambda_{H}|+\|u-u_{H}\|_{-\frac{1}{2},\partial\Omega}.
\end{eqnarray}
Since
\begin{eqnarray*}
\|P_{h}u-u^{w,h}\|_{0,\Omega_{0}}&\leq & \|P_{h}u-u^{w}\|_{0,\Omega_{0}}+\|e^{h}\|_{0,\Omega_{0}}\\
&\leq &\|P_{h}u-P_{w}u\|_{0,\Omega_{0}}+\|P_{w}u-u^{w}\|_{0,\Omega_{0}}+\|e^{h}\|_{0,\Omega_{0}},
\end{eqnarray*}
which together with (\ref{s5.4}) and  (\ref{s5.6}) we get
\begin{eqnarray}\label{s5.7}
\|P_{h}u-u^{w,h}\|_{1,D}&\lesssim& |\lambda-\lambda_{H}|+\|u-u_{H}\|_{-\frac{1}{2},\partial\Omega}+K\|P_{w}u-u^{w}\|_{0,\Omega}\nonumber\\
&&+\|P_{h}u-P_{w}u\|_{0,\Omega}+\|e^{h}\|_{0,\Omega_{0}}.
\end{eqnarray}

Next, we will use Aubin-Nitsche duality argument to estimate $\|e^{h}\|_{0,\Omega_{0}}$. Given any $\zeta\in L^{2}(\Omega_{0})$,
there exists $\eta\in H^{1}_{\Gamma}(\Omega_{0})$ satisfying (\ref{s3.2}), namely,
$$a(v,\eta)=(v,\zeta),~~~\forall v\in H^{1}_{\Gamma}(\Omega_{0}).$$
Let $\eta_{h}^{0}\in V_{h}^{0}(\Omega_{0})$ and
$\eta_{H}^{0}\in V_{H}^{0}(\Omega_{0})$ satisfy
\begin{eqnarray*}
a(v,\eta-\eta_{h}^{0})=0,~~\forall v\in V_{h}^{0}(\Omega_{0}),\\
a(v,\eta-\eta_{H}^{0})=0,~~\forall v\in V_{H}^{0}(\Omega_{0}).
\end{eqnarray*}
Then we deduce that
\begin{eqnarray*}
(e^{h}, \zeta)&=&a(e^{h},\eta)=a(e^{h},\eta^{0}_{h})=a(u^{w,h}-u^{w},\eta^{0}_{h})\\
&=&a(P_{h}u-u^{w},\eta^{0}_{h})+a(u^{w,h},\eta^{0}_{h})-a(P_{h}u,\eta^{0}_{h})\\
&=&a(P_{h}u-u^{w},\eta^{0}_{h})-\lambda_{H}\langle u_{H}, \eta_{h}^{0}\rangle+\lambda \langle u,\eta_{h}^{0}\rangle\\
&=&a(P_{h}u-u^{w},\eta^{0}_{h}-\eta)+a(P_{h}u-u^{w},\eta-\eta^{0}_{H})+a(P_{h}u-u^{w},\eta^{0}_{H})\\
&&~~~-\lambda_{H}\langle u_{H}, \eta_{h}^{0}\rangle+\lambda \langle u,\eta_{h}^{0}\rangle\\
&=&a(P_{h}u-u^{w},\eta^{0}_{h}-\eta)+a(P_{h}u-u^{w},\eta-\eta^{0}_{H})\\
&&~~~-\lambda\langle u, \eta_{H}^{0}\rangle+\lambda_{H} \langle u_{H},\eta_{H}^{0}\rangle-\lambda_{H}\langle u_{H}, \eta_{h}^{0}\rangle+\lambda \langle u,\eta_{h}^{0}\rangle\\
&=&a(P_{h}u-u^{w},\eta^{0}_{h}-\eta)+a(P_{h}u-u^{w},\eta-\eta^{0}_{H})\\
&&~~~+\lambda\langle u, \eta_{h}^{0}-\eta_{H}^{0}\rangle+\lambda_{H} \langle u_{H},\eta_{H}^{0}-\eta_{h}^{0}\rangle\\
&=&a(P_{h}u-u^{w},\eta^{0}_{h}-\eta)+a(P_{h}u-u^{w},\eta-\eta^{0}_{H})\\
&&~~~+\langle \lambda u-\lambda_{H} u_{H}, \eta_{h}^{0}-\eta_{H}^{0}\rangle.
\end{eqnarray*}
From the error estimate of finite element and the local regularity assumption {\bf $R(\Omega_{0})$} we have
\begin{eqnarray*}
\|\eta^{0}_{H}-\eta\|_{1,\Omega_{0}}\lesssim H^{r}\|\zeta\|_{0,\Omega_{0}},~~~\|\eta^{0}_{h}-\eta\|_{1,\Omega_{0}}\lesssim h^{r}\|\zeta\|_{0,\Omega_{0}}.
\end{eqnarray*}
Thus, we get the estimation for any $\zeta\in L^{2}(\Omega_{0})$ that
\begin{eqnarray*}
|(e^{h}, \zeta)|\leq (H^{r}\|P_{h}u-u^{w}\|_{1,\Omega}+|\lambda-\lambda_{H}|+\|u-u_{H}\|_{-\frac{1}{2},\partial\Omega})\|\zeta\|_{0,\Omega_{0}},
\end{eqnarray*}
which leads to
\begin{eqnarray*}
\|e^{h}\|_{0,\Omega_{0}}\lesssim |\lambda-\lambda_{H}|+\|u-u_{H}\|_{-\frac{1}{2},\partial\Omega}+H^{r}\|P_{h}u-u^{w}\|_{1,\Omega}.
\end{eqnarray*}
From (\ref{s5.4}) and the triangle inequality
\begin{eqnarray*}
\|P_{h}u-u^{w}\|_{1,\Omega}\leq \|P_{h}u-P_{w}u\|_{1,\Omega}+\|P_{w}u-u^{w}\|_{1,\Omega},
\end{eqnarray*}
we obtain
\begin{eqnarray}\label{s5.8}
\|e^{h}\|_{0,\Omega_{0}}&\lesssim &|\lambda-\lambda_{H}|+\|u-u_{H}\|_{-\frac{1}{2},\partial\Omega}+H^{r}\|P_{h}u-P_{w}u\|_{1,\Omega}\nonumber\\
&&~~+K H^{r}\|P_{w}u-u^{w}\|_{0,\Omega}.
\end{eqnarray}
From Aubin-Nitsche duality argument we can easily get the estimation
\begin{eqnarray*}
\|P_{h}u-P_{w}u\|_{0,\Omega}\lesssim w^{r}\|P_{h}u-P_{w}u\|_{1,\Omega}.
\end{eqnarray*}
Substituting the above estimate and (\ref{s5.8}) into (\ref{s5.7}) we obtain
\begin{eqnarray*}
\|P_{h}u-u^{w,h}\|_{1,D}&\lesssim & |\lambda-\lambda_{H}|+\|u-u_{H}\|_{-\frac{1}{2},\partial\Omega}+H^{r}\|P_{h}u-P_{w}u\|_{1,\Omega}\\
&~~~&+w^{r}\|P_{h}u-P_{w}u\|_{1,\Omega}+K(H^{r}+1)\|P_{w}u-u^{w}\|_{0,\Omega}\\
&\lesssim & |\lambda-\lambda_{H}|+\|u-u_{H}\|_{-\frac{1}{2},\partial\Omega}+H^{r}\|P_{h}u-P_{w}u\|_{1,\Omega}\\
&~~~&+\|P_{w}u-u^{w}\|_{0,\Omega},
\end{eqnarray*}
which together with Theorem 2.1, Lemma 2.2, Lemma 2.4 and Theorem 4.1 yields
\begin{eqnarray}\label{s5.9}
\|P_{h}u-u^{w,h}\|_{1,D}&\lesssim & H^{\frac{2r}{\alpha}}+H^{\frac{r+r}{\alpha}}+H^{r}w^{r}+w^{r+r}+H^{\frac{r+r}{\alpha}}\lesssim H^{\frac{2r}{\alpha}}.
\end{eqnarray}
\indent Similarly, by using the same argument for $(G\setminus \overline{D})\subset\subset \Omega_{0}$ we can obtain the estimation
\begin{eqnarray}\label{s5.10}
\|P_{h}u-u^{w,h}\|_{1,G\setminus \overline{D}}\lesssim  H^{\frac{2r}{\alpha}}.
\end{eqnarray}

\indent Now, the remainder is to analyze $\|P_{h}u-u^{w,h}\|_{1,\Omega\setminus \overline{G}}$. From the definition of $u^{w,h}$ we see that
\begin{eqnarray*}
\|P_{h}u-u^{w,h}\|_{1,\Omega\setminus \overline{\Omega_{0}}}=\|P_{h}u-u^{w}\|_{1,\Omega\setminus \overline{\Omega_{0}}},
\end{eqnarray*}
thus
\begin{eqnarray*}
\|P_{h}u-u^{w,h}\|_{1,\Omega\setminus \overline{G}}&\leq& \|P_{h}u-u^{w}\|_{1,\Omega\setminus \overline{\Omega_{0}}}+\|P_{h}u-u^{w}\|_{1,\Omega_{0}\setminus \overline{G}}+\|e^{h}\|_{1,\Omega_{0}\setminus \overline{G}}\\
&\lesssim & \|P_{h}u-u^{w}\|_{1,\Omega\setminus \overline{G}}+\|e^{h}\|_{1,\Omega_{0}\setminus \overline{G}}\\
&\lesssim &\|P_{h}u-u\|_{1,\Omega\setminus \overline{G}}+\|u-u^{w}\|_{1,\Omega\setminus \overline{G}}
+\|e^{h}\|_{1,\Omega_{0}\setminus \overline{G}}.
\end{eqnarray*}
From the fact that
\begin{eqnarray}\label{s5.11}
a(e^{h},v)=-\lambda_{H}\langle u_{H},v\rangle+\lambda \langle u,v\rangle-a(u^{w}-u,v),~~~\forall v\in V_{0}^{h}(\Omega_{0}\setminus \overline{F})
\end{eqnarray}
and Lemma 3.2 we derive that
\begin{eqnarray}\label{s5.12}
\|e^{h}\|_{1,\Omega_{0}\setminus \overline{G}}\lesssim \|e^{h}\|_{0,\Omega_{0}\setminus \overline{F}}+|\lambda-\lambda_{H}|+\|u-u_{H}\|_{-\frac{1}{2},\partial\Omega}
+\|u^{w}-u\|_{1,\Omega_{0}\setminus \overline{F}};
\end{eqnarray}
we then arrive at
\begin{eqnarray}\label{s5.13}
\|P_{h}u-u^{w,h}\|_{1,\Omega\setminus \overline{G}}&\lesssim& \|e^{h}\|_{0,\Omega_{0}}+|\lambda-\lambda_{H}|+\|u-u_{H}\|_{-\frac{1}{2},\partial\Omega}\nonumber\\
&~&+\|P_{h}u-u\|_{1,\Omega\setminus \overline{F}}+\|u^{w}-u\|_{1,\Omega\setminus \overline{F}}.
\end{eqnarray}
From (\ref{s5.4}) and the triangle inequality
\begin{eqnarray*}
\|u^{w}-u\|_{1,\Omega\setminus \overline{F}}\leq \|u-P_{w}u\|_{1,\Omega\setminus \overline{F}}+\|P_{w}u-u^{w}\|_{1,\Omega\setminus \overline{F}},
\end{eqnarray*}
we get
\begin{eqnarray*}
\|u^{w}-u\|_{1,\Omega\setminus \overline{F}}&\leq &\|u-P_{w}u\|_{1,\Omega\setminus \overline{F}}+|\lambda-\lambda_{H}|+\|u-u_{H}\|_{-\frac{1}{2},\partial\Omega}\\
&&~~+K\|P_{w}u-u^{w}\|_{0,\Omega}.
\end{eqnarray*}
Thus, substituting the above inequality and (\ref{s5.8}) into (\ref{s5.13}) we conclude that
\begin{eqnarray}\label{s5.14}
&&\|P_{h}u-u^{w,h}\|_{1,\Omega\setminus \overline{G}}\lesssim |\lambda-\lambda_{H}|+\|u-u_{H}\|_{-\frac{1}{2},\partial\Omega}+H^{r}\|P_{h}u-P_{w}u\|_{1,\Omega}\nonumber\\
&&~~~~~~+K (H^{r}+1)\|P_{w}u-u^{w}\|_{0,\Omega}+\|u-P_{h}u\|_{1,\Omega\setminus \overline{F}}+\|u-P_{w}u\|_{1,\Omega\setminus \overline{F}}.
\end{eqnarray}
Since $(\Omega\setminus \overline{F})\subset\subset (\Omega\setminus \overline{D})$, we obtain by using Theorem 3.2 that
\begin{eqnarray*}
\|P_{h}u-u^{w,h}\|_{1,\Omega\setminus \overline{G}}&\lesssim& |\lambda-\lambda_{H}|+\|u-u_{H}\|_{-\frac{1}{2},\partial\Omega}+H^{r}\|P_{h}u-P_{w}u\|_{1,\Omega}\nonumber\\
&~&+K (H^{r}+1)\|P_{w}u-u^{w}\|_{0,\Omega}\\
&~&+\inf\limits_{v\in V_{h}(\Omega)}\|u-v\|_{1,\Omega\setminus \overline{D}}+\|u-P_{h}u\|_{0,\Omega}\nonumber\\
&~&+\inf\limits_{v\in V_{w}(\Omega)}\|u-v\|_{1,\Omega\setminus \overline{D}}+\|u-P_{w}u\|_{0,\Omega},
\end{eqnarray*}
which together with Theorem 2.1, Lemma 2.2, Lemma 2.4 and Theorem 4.1 yields
\begin{eqnarray}\nonumber
&&\|P_{h}u-u^{w,h}\|_{1,\Omega\setminus \overline{G}}\lesssim H^{\frac{2r}{\alpha}}+H^{\frac{r+r}{\alpha}}+H^{r}w^{r}+w^{r}w^{r}+H^{\frac{r+r}{\alpha}}+w+w^{r}w^{r}\\\label{s5.15}
&&~~~\lesssim H^{\frac{2r}{\alpha}}+w.
\end{eqnarray}
Combining (\ref{s5.15}), (\ref{s5.10}), (\ref{s5.9}) and (\ref{s2.20}), we prove the desired result (\ref{s5.1}).\\
\indent Similarly, we can prove (\ref{s5.1a}). Using the same argument of (\ref{s4.4a}) we can prove (\ref{s5.2}).~~~$\Box$\\

\indent In the first two steps of Scheme 3  we actually use Scheme 1 to compute. We can also use the second kind of two-grid scheme, Scheme 2, to devise the local computational scheme. When the number of isolated singular points is larger than 1, we can design the parallel version of Scheme 3.\\

\section{Numerical experiments}
In this section, we shall report some numerical experiments to show the efficiency of our schemes.
Consider the problem (\ref{s2.1}) with $k=1$ on the test domain $\Omega_{S}=(-\frac{\sqrt{2}}{2}, \frac{\sqrt{2}}{2})^2$, $\Omega_{L}=(-1, 1)^2\setminus([0, 1)\times (-1, 0])$, and $\Omega_{SL}=(-\frac{\sqrt{2}}{2}, \frac{\sqrt{2}}{2})^2\setminus\{0\leq x\leq\frac{\sqrt{2}}{2}, y=0\}$.
We use Matlab 2012a to solve (\ref{s2.1}) on a Lenovo ideaPad PC with 1.8GHZ CPU and 8GB RAM. Our program is compiled under the
package of iFEM \cite{chen}.  In our computation, we adopt a uniform isosceles right triangulation and the triangle linear element ($m=1$), and take $n(x)=4$ or $n(x)=4+4i$. For simplicity we use the following notations in our tables :\\
\indent $S,L, and ~Slit$ stand for the domain $\Omega_{S}, \Omega_{L}$ and $\Omega_{SL}$, respectively.\\
\indent$\lambda_{j,H}$ is the $j$th approximate eigenvalue derived from Step 1 of our schemes. Here we use the sparse solver
$eigs(A,B,j,'sm')$ to get the first $j$ eigenvalues.\\
\indent $\lambda_{j}^{w}$ is the $j$th approximate eigenvalue obtained by Step 2 of Schemes 1 and 2. \\
Here we use Matlab solver $'\setminus'$ to solve $j$ equations at the same time to get the first $j$ eigenvalues. \\
\indent $\lambda_{j}^{w,h}$ is the $j$th approximate eigenvalue  obtained by Scheme 3.\\
\indent $dof_H$ is the number of degrees of freedom for solving the eigenvalue problem directly on $\pi_H(\Omega)$.\\
\indent $dof_w$ is the number of degrees of freedom for solving the boundary problem on mesoscopic grid $\pi_w(\Omega)$.\\
\indent $dof_h$ is the number of degrees of freedom for solving the boundary problem on locally fine mesh $\pi_h(\Omega_0)$.\\
\indent $t$(s) is the CPU time from the program starting to the current calculating results appearing.\\
\indent The symbol '--' means that the calculation cannot proceed since the computer runs out of memory.\\

\indent  According to the regularity results, we have
$r=1$ on $\Omega_{S}$, $r=\frac{2}{3}$ on $\Omega_{L}$, and $r=\frac{1}{2}$ on $\Omega_{SL}$. Thus, the approximate eigenvalues obtained by Schemes 1 and 2, when taking $w=\mathcal {O}(H^2)$, can achieve $\mathcal {O}(w^2)$ on $\Omega_{S}$, $\mathcal {O}(w^{\frac{4}{3}})$ on $\Omega_{L}$ and $\mathcal {O}(w)$ on $\Omega_{SL}$.
When the index of refraction $n(x)$ is real, the problem is selfadjoint and all Stekloff eigenvalues are real. Comparing Table 1 and Tables 2-3, Table 4 and Tables 5 and 6 we can see that under the same accuracy, the two-grid discretization Schemes 1 and 2 take less time to get the asymptotically optimal approximations. Especially, Scheme 2 works more efficiently than directly solving and Scheme 1 since
the matrices are constructed to be symmetric and definite in solving linear systems.
\\
From numerical experiments we find that the eigenfunction corresponding to the second eigenvalue is singular near the origin on $\Omega_{L}$ and $\Omega_{SL}$. We compute the second approximate eigenvalue by Scheme 3 on $\Omega_{L}$ and $\Omega_{SL}$, and the results are listed in Tables 7  and 8 from which we can see that the local correction does work.\\

\begin{table}[htbp]
  \centering
  \caption{The first four eigenvalues computed directly when $n=4$.}
    \begin{tabular}{ccccccc}
    \hline\noalign{\smallskip}
    ~domain~&$H$ &$\lambda_{1,H}$ & $\lambda_{2,H}$ & $\lambda_{3,H}$ & $\lambda_{4,H}$&$t$ \\
    \noalign{\smallskip}\hline\noalign{\smallskip}
    S     & 1/512 & 2.202501387  & -0.212254531  & -0.212255107  & -0.908066632  & 20.77  \\
    S     & 1/1024 & 2.202505691  & -0.212252760  & -0.212252904  & -0.908058722  & 109.80  \\
    L     & 1/512 & 2.533187700  & 0.857690917  & 0.124518848  & -1.085315271  & 14.14  \\
    L     & 1/1024 & 2.533207148  & 0.857750492  & 0.124523033  & -1.085302932  & 67.98  \\
    Slit  & 1/512 & 1.484704242  & 0.460698784  & -0.184178326  & -0.690081852  & 20.46  \\
    Slit  & 1/1024 & 1.484709990  & 0.461215008  & -0.184176518  & -0.690076859  & 98.87  \\
 \hline
    \end{tabular}%
\end{table}%

\begin{table}[htbp]
  \centering
  \caption{The first four eigenvalues computed by Scheme 1 when $n=4$.}
    \begin{tabular}{cccccccc}
    \hline\noalign{\smallskip}
    ~domain~&$H$ &$w$ &$\lambda_{1}^{w}$ & $\lambda_{2}^{w}$ & $\lambda_{3}^{w}$ & $\lambda_{4}^{w}$&$t$ \\
    \noalign{\smallskip}\hline\noalign{\smallskip}
S     &   1/64 & 1/512 & 2.202501132  & -0.212254531  & -0.212255108  & -0.908066630  & 15.71  \\
    S     &   1/64 & 1/1024 & 2.202505431  & -0.212252760  & -0.212252904  & -0.908058720  & 82.38  \\
    S     &    1/128 & 1/1024 & 2.202505676  & -0.212252760  & -0.212252904  & -0.908058722  & 83.45  \\
    L     &   1/64 & 1/512 & 2.533179767  & 0.857690001  & 0.124518848  & -1.085313547  & 10.22  \\
    L     &   1/64 & 1/1024 & 2.533199015  & 0.857749498  & 0.124523032  & -1.085301167  & 49.58  \\
    L     &    1/128 & 1/1024 & 2.533206625  & 0.857750348  & 0.124523033  & -1.085302824  & 50.07  \\
    Slit  &   1/64 & 1/512 & 1.484704002  & 0.460697127  & -0.184178327  & -0.690081762  & 11.91  \\
    Slit  &   1/64 & 1/1024 & 1.484709743  & 0.461213101  & -0.184176518  & -0.690076761  & 58.80  \\
    Slit  &    1/128 & 1/1024 & 1.484709973  & 0.461214611  & -0.184176518  & -0.690076844  & 59.84  \\
 \hline
    \end{tabular}%
\end{table}%

\begin{table}[htbp]
  \centering
  \caption{The first four eigenvalues computed by Scheme 2 when $n=4$.}
    \begin{tabular}{cccccccc}
    \hline\noalign{\smallskip}
    ~domain~&$H$ &$w$ &$\lambda_{1}^{w}$ & $\lambda_{2}^{w}$ & $\lambda_{3}^{w}$ & $\lambda_{4}^{w}$&$t$ \\
    \noalign{\smallskip}\hline\noalign{\smallskip}
S     & 1/64  & 1/512 & 2.2025013  & -0.2122545  & -0.2122551  & -0.9080666  & 8.54  \\
    S     & 1/64  & 1/1024 & 2.2025056  & -0.2122528  & -0.2122529  & -0.9080587  & 41.73  \\
    S     & 1/128 & 1/1024 & 2.2025057  & -0.2122528  & -0.2122529  & -0.9080587  & 42.60  \\
    L     & 1/64  & 1/512 & 2.5331872  & 0.8576891  & 0.1245188  & -1.0853154  & 6.11  \\
    L     & 1/64  & 1/1024 & 2.5332066  & 0.8577485  & 0.1245230  & -1.0853030  & 27.53  \\
    L     & 1/128 & 1/1024 & 2.5332071  & 0.8577502  & 0.1245230  & -1.0853029  & 28.08  \\
    Slit  & 1/64  & 1/512 & 1.4847042  & 0.4606756  & -0.1841783  & -0.6900818  & 8.21  \\
    Slit  & 1/64  & 1/1024 & 1.4847099  & 0.4611884  & -0.1841765  & -0.6900768  & 37.75  \\
    Slit  & 1/128 & 1/1024 & 1.4847100  & 0.4612092  & -0.1841765  & -0.6900768  & 38.82  \\
 \hline
    \end{tabular}%
\end{table}%

\begin{table}[htbp]
  \centering
  \caption{The first four eigenvalues computed directly when $n=4+4i$.}
    \begin{tabular}{ccccccc}
    \hline\noalign{\smallskip}
    ~domain~&$H$ &$\lambda_{1,H}$ & $\lambda_{2,H}$ & $\lambda_{3,H}$ & $\lambda_{4,H}$&$t$ \\
    \noalign{\smallskip}\hline\noalign{\smallskip}
    S     & 1/64  & 0.686951  & -0.343131  & -0.342924  & -2.802148  & 0.75  \\
          &       & +2.495332i & +0.850617i & +0.85054i & +0.542231i &  \\
    S     & 1/128 & 0.686652  & -0.343068  & -0.343016  & -2.797931  & 1.81  \\
          &       & +2.495304i & +0.850714i &  +0.850695i&  +0.541106i&  \\
    S     & 1/256 & 0.686577  & -0.343052  & -0.343039  & -2.796876  & 9.42  \\
          &       &  +2.495296i&  +0.850738i&  +0.850734i&  +0.540824i&  \\
    S     & 1/512 & 0.686558  & -0.343048  & -0.343045  & -2.796612  & 46.18  \\
          &       &  +2.495295i&  +0.850744i&  +0.850743i&  +0.540753i&  \\
    S     & 1/1024 &   --    & --      &  --     &  --     & -- \\
    L     & 1/64  & 0.5163544 & 0.39617526 & -0.0769975 & -1.4419097 & 0.60  \\
          &       &  +2.882867i&  +1.457866i&  +1.04222i&  +0.805745i&  \\
    L     & 1/128 & 0.5148057 & 0.39665783 & -0.0771338 & -1.4408607 & 1.27  \\
          &       &  +2.882465i&  +1.458552i&  +1.042563i&  +0.804959i&  \\
    L     & 1/256 & 0.5144169 & 0.39687654 & -0.0771676 & -1.4405978 & 6.98  \\
          &       &  +2.882359i&  +1.458814i&  +1.042649i&  +0.804761i&  \\
    L     & 1/512 & 0.5143195 & 0.39696985 & -0.077176 & -1.4405319 & 33.05  \\
          &       &  +2.882332i&  +1.458916i&  +1.042671i&  +0.804711i&  \\
    L     & 1/1024 & --      &  --     &  --     & --      & -- \\
    Slit  & 1/64  & 0.9198804 & 0.28552179 & -0.2626473 & -0.7423903 & 0.96  \\
          &       &  +1.770436i&  +0.995916i&  +0.75731i&  +0.608702i&  \\
    Slit  & 1/128 & 0.9194638 & 0.28906713 & -0.2626227 & -0.7421686 & 1.87  \\
          &       &  +1.770697i&  +0.997926i&  +0.757415i&  +0.608759i&  \\
    Slit  & 1/256 & 0.9193482 & 0.29084621 & -0.2626166 & -0.7421124 & 9.53  \\
          &       &  +1.770765i&  +0.998908i&  +0.757442i&  +0.608773i&  \\
    Slit  & 1/512 & 0.9193164 & 0.29173723 & -0.2626151 & -0.7420981 & 46.20  \\
          &       &  +1.770782i&  +0.999395i&  +0.757448i&  +0.608775i&  \\
    Slit  & 1/1024 &  --     &   --    &  --     &  --     & -- \\
 \hline
    \end{tabular}%
\end{table}%

\begin{table}[htbp]
  \centering
  \caption{The first four eigenvalues computed by Scheme 1 when $n=4+4i$.}
    \begin{tabular}{cccccccc}
    \hline\noalign{\smallskip}
    ~domain~&$H$&$w$ &$\lambda_{1}^{w}$ & $\lambda_{2}^{w}$ & $\lambda_{3}^{w}$ & $\lambda_{4}^{w}$&$t$ \\
    \noalign{\smallskip}\hline\noalign{\smallskip}
 S &   1/64 &  1/512& 0.6865577  & -0.3430479  & -0.3430446  & -2.7966123  & 39.95  \\
          &       &       & +2.4952946i & +0.8507445i & +0.8507433i & +0.5407542i &  \\
     S&   1/64 &  1/1024& 0.6865530  & -0.3430469  & -0.3430461  & -2.7965463  & 230.36  \\
          &       &       & +2.4952942i & +0.850746i & +0.8507457i & +0.5407366i &  \\
     S&    1/128 &  1/1024& 0.6865534  & -0.3430469  & -0.3430461  & -2.7965461  & 215.01  \\
          &       &       & +2.4952941i & +0.850746i & +0.8507457i & +0.5407357i &  \\
     L&   1/64 &  1/512& 0.5143181  & 0.3969728  & -0.0771760  & -1.4405317  & 25.28  \\
          &       &       & +2.8823326i & +1.4589166i & +1.0426708i & +0.8047119i &  \\
     L&   1/64 &  1/1024& 0.5142937  & 0.3970116  & -0.0771781  & -1.4405152  & 131.40  \\
          &       &       & +2.8823258i & +1.4589567i & +1.0426763i & +0.8046995i &  \\
     L&    1/128 &  1/1024& 0.5142950  & 0.3970089  & -0.0771780  & -1.4405153  & 133.61  \\
          &       &       & +2.8823255i & +1.4589562i & +1.0426762i & +0.804699i &  \\
     Slit&   1/64 &  1/512& 0.9193164  & 0.2917455  & -0.2626151  & -0.7420981  & 34.47  \\
          &       &       & +1.7707824i & +0.9993949i & +0.7574481i & +0.6087755i &  \\
     Slit&   1/64 &  1/1024& 0.9193077  & 0.2921926  & -0.2626147  & -0.7420944  & 181.84  \\
          &       &       & +1.770787i & +0.999637i & +0.7574498i & +0.608776i &  \\
     Slit&    1/128 &  1/1024& 0.9193078  & 0.2921851  & -0.2626147  & -0.7420944  & 187.01  \\
          &       &       & +1.770787i & +0.999637i & +0.7574498i & +0.608776i &  \\
 \hline
    \end{tabular}%
\end{table}%

\begin{table}[htbp]
  \centering
  \caption{The first four eigenvalues computed by Scheme 2 when $n=4+4i$.}
    \begin{tabular}{cccccccc}
    \hline\noalign{\smallskip}
    ~domain~&$H$&$w$ &$\lambda_{1}^{w}$ & $\lambda_{2}^{w}$ &$\lambda_{3}^{w}$ & $\lambda_{4}^{w}$&$t$ \\
    \noalign{\smallskip}\hline\noalign{\smallskip}
    S &  1/64&  1/512& 0.6866414 & -0.3430158 & -0.3429982 & -2.7957018 & 14.19  \\
          &       &       & +2.4955259i & +0.8507372i & +0.8507316i & +0.5401331i &  \\
     S&  1/64&  1/1024&  0.6865532 &  -0.3430469 &  -0.3430460 &  -2.7965460 & 67.51  \\
          &       &       & +2.4955282i & +0.8507386i & +0.8507338i & +0.5401081i &  \\
     S&  1/128&  1/1024&  0.6865743 &  -0.3430388 &  -0.3430345 &  -2.7963187 & 69.35  \\
          &       &       & +2.4953521i & +0.8507442i & +0.8507428i & +0.54058i &  \\
     L&  1/64&  1/512&  0.5148525 &  0.3975831 &  -0.0770520 &  -1.4407879 & 9.97  \\
          &       &       & +2.8832036i & +1.4581667i & +1.0426418i & +0.8041698i &  \\
     L&  1/64&  1/1024&  0.5148345 &  0.3976471 &  -0.0770526 &  -1.4407746 & 44.34  \\
          &       &       & +2.8832075i & +1.458175i & +1.0426468i & +0.8041508i &  \\
     L&  1/128&  1/1024&  0.5144287 &  0.3972538 &  -0.0771471 &  -1.4405817 & 46.03  \\
          &       &       & +2.8825473i & +1.4586488i & +1.0426688i & +0.8045621i &  \\
     Slit&  1/64&  1/512&  0.9194996 &  0.2930070 &  -0.2625644 &  -0.7420884 & 14.15  \\
          &       &       & +1.7708477i & +0.9969382i & +0.7574283i & +0.6087121i &  \\
     Slit&  1/64&  1/1024&  0.9194939 &  0.2935419 &  -0.2625634 &  -0.7420857 & 62.20  \\
          &       &       & +1.7708533i & +0.9970054i & +0.7574298i & +0.6087082i &  \\
     Slit&  1/128&  1/1024&  0.9193597 &  0.2928178 &  -0.2626020 &  -0.7420970 & 63.15  \\
          &       &       & +1.7708054i & +0.9984062i & +0.7574448i & +0.6087449i &  \\
 \hline
    \end{tabular}%
\end{table}%

\begin{table}[htbp]
  \centering
  \caption{The second eigenvalue computed by Scheme 3 when $n=4$.}
    \begin{tabular}{cccccccc}
    \hline\noalign{\smallskip}
    ~domain~&$dof_H$ &$dof_w$ &$dof_h$& $\lambda_{2,H}$ & $\lambda_{2}^{w}$ & $\lambda_{2}^{w,h}$ &$t$ \\
    \noalign{\smallskip}\hline\noalign{\smallskip}
     L     & 3201  & 49665 & 48896 & 0.8561269  & 0.8575382  & 0.8576781  & 5.60  \\
    L     & 3201  & 49665 & 196096 & 0.8561269  & 0.8575382  & 0.8577346  & 11.32  \\
    L     & 3201  & 49665 & 785408 & 0.8561269  & 0.8575382  & 0.8577574  & 42.08  \\
    Slit  & 4257  & 66177 & 65152 & 0.4533833  & 0.4596629  & 0.4606905  & 7.15  \\
    Slit  & 4257  & 66177 & 261376 & 0.4533833  & 0.4596629  & 0.4612043  & 14.92  \\
    Slit  & 4257  & 66177 & 1047040 & 0.4533833  & 0.4596629  & 0.4614613  & 53.20  \\
 \hline
    \end{tabular}%
\end{table}%

\begin{table}[htbp]
  \centering
  \caption{The second eigenvalue computed by Scheme 3 when $n=4+4i$.}
    \begin{tabular}{cccccccc}
    \hline\noalign{\smallskip}
    ~domain~&$dof_H$ &$dof_w$ &$dof_h$& $\lambda_{2,H}$ & $\lambda_{2}^{w}$ & $\lambda_{2}^{w,h}$ &$t$ \\
    \noalign{\smallskip}\hline\noalign{\smallskip}
    L& 3201  & 49665 & 48896 & 0.3961753 & 0.3968789 & 0.3969791 & 7.95  \\
          &       &       &       & +1.4578657i & +1.4588146i & +1.458904i &  \\
    L& 3201  & 49665 & 196096 & 0.3961753 & 0.3968789 & 0.3970194 & 23.14  \\
          &       &       &       & +1.4578657i & +1.4588146i & +1.4589409i &  \\
    L& 3201  & 49665 & 785408 & 0.3961753 & 0.3968789 & 0.3970356 & 115.94  \\
          &       &       &       & +1.4578657i & +1.4588146i & +1.458956i &  \\
    Slit& 4257  & 66177 & 65152 & 0.2855218 & 0.2908523 & 0.2917435 & 11.08  \\
          &       &       &       & +0.9959158i & +0.9989086i & +0.9993878i &  \\
    Slit& 4257  & 66177 & 261376 & 0.2855218 & 0.2908523 & 0.2921890 & 32.71  \\
          &       &       &       & +0.9959158i & +0.9989086i & +0.9996281i &  \\
    Slit& 4257  & 66177 & 1047040 & 0.28552179& 0.29085229& 0.2924118& 175.7167 \\
          &       &       &       & +0.9959158i & +0.9989086i & +0.9997486i &  \\
 \hline
    \end{tabular}%
\end{table}%



\noindent {\bf Acknowledgements} This work is supported by National
Natural Science Foundation of China(Grant No.11761022).



\end{document}